\documentclass[microtype]{gtpart}   
\usepackage{latexsym, amsfonts, amscd, amsmath, verbatim, graphicx, amsthm}

\title[Floer homology of surface symplectomorphisms]{Symplectic Floer homology of area-preserving surface diffeomorphisms}

\author[A. Cotton-Clay]{Andrew Cotton-Clay}
\givenname{Andrew}
\surname{Cotton-Clay}
\address{Department of Mathematics\\University of California, Berkeley\\\newline
         Berkeley, CA 94720-3840\\USA}
\email{acotton@math.berkeley.edu}
\urladdr{http://math.berkeley.edu/~acotton}

\keyword{Floer homology}
\keyword{symplectomorphism}
\keyword{surface diffeomorphism}
\keyword{mapping class group}
\keyword{fixed point}
\keyword{Nielsen class}
\subject{primary}{msc2000}{53D40}
\subject{primary}{msc2000}{37J10}

\arxivreference{arXiv:0807.2488}  
\arxivpassword{stpz3}   

\def\be{\begin{enumerate}}
\def\ee{\end{enumerate}}
\def\ii{\item}
\def\bi{\begin{itemize}}
\def\ei{\end{itemize}}

\def\RR{\mathbb{R}}
\def\CC{\mathbb{C}}

\def\beqa{\begin{eqnarray*}}
\def\eeqa{\end{eqnarray*}}
\def\beqan{\begin{eqnarray}}
\def\eeqan{\end{eqnarray}}

\def\ZZ{\mathbb{Z}}
\def\PP{\mathbb{P}}

\def\MM{\mathcal{M}}
\def\TT{\mathcal{T}}

\def\HH{\mathbb{H}}

\def\del{\partial}

\newtheorem{thm}{Theorem}[section]
\newtheorem{cond}[thm]{Condition}
\newtheorem{lemma}[thm]{Lemma}

\newtheorem{prop}[thm]{Proposition}
\newtheorem{cor}[thm]{Corollary}

\theoremstyle{definition}
\newtheorem{defn}[thm]{Definition}

\newtheorem{remark}[thm]{Remark}

\def\proof{\textit{Proof:}\ \ \ \ }

\def\eop{\rule{5pt}{5pt}}
\def\qed{\hfill $\eop$ \vskip10pt}

\def\Hom{\textrm{Hom}}

\def\Re{\textrm{Re}}
\def\Im{\textrm{Im}}

\def\det{\textrm{det}}

\def\lim{\textrm{lim}}

\def\dbar{\overline\partial}

\def\ker{\textrm{ker}}
\def\ind{\textrm{ind}}

\def\id{\textrm{id}}

\def\Symp{\textrm{Symp}}
\def\Diff{\textrm{Diff}}
\def\tr{\textrm{tr}}

\begin{document}

\begin{abstract}
The symplectic Floer homology $HF_*(\phi)$ of a symplectomorphism $\phi: \Sigma \rightarrow \Sigma$ encodes data about the fixed points of $\phi$ using counts of holomorphic cylinders in $\RR \times M_\phi$, where $M_\phi$ is the mapping torus of $\phi$. We give an algorithm to compute $HF_*(\phi)$ for $\phi$ a surface symplectomorphism in a pseudo-Anosov or reducible mapping class, completing the computation of Seidel's $HF_*(h)$ for $h$ any orientation-preserving mapping class.
\end{abstract}

\begin{asciiabstract}
The symplectic Floer homology HF_*(f) of a symplectomorphism f: S -> S encodes data about the fixed points of f using counts of holomorphic cylinders in R x M_f, where M_f is the mapping torus of f. We give an algorithm to compute HF_*(f) for f a surface symplectomorphism in a pseudo-Anosov or reducible mapping class, completing the computation of Seidel's HF_*(h) for h any orientation-preserving mapping class.
\end{asciiabstract}

\maketitle

\section{Introduction}

\subsection{Statement of results}
\label{results}

Let $\Sigma$ be a compact, connected, oriented surface.
Let $\omega$ be a symplectic form (here, an area form) on $\Sigma$. Let $\phi \in \Symp(\Sigma,\omega)$ be a symplectomorphism (here, an area-preserving diffeomorphism) of $\Sigma$ with nondegenerate fixed points. That is, the fixed points of $\phi$ are cut out transversally in the sense that $\det(1-d\phi_x) \neq 0$ for fixed points $x$. We then consider the symplectic Floer homology chain complex $$CF_*(\phi) := \ZZ/2 \left<x \left| x \in \textrm{Fix}(\phi)\right.\right>,$$ a $\ZZ/2$-graded chain complex over $\ZZ/2$. The grading of a generator $x$ is given by the sign of $\det(1-d\phi_x)$. Note that a fixed point $x$ of $\phi$ corresponds to a constant section of the mapping torus $M_\phi$, the $\Sigma$-bundle over $S^1$ with monodromy $\phi$. The matrix coefficient $\left<\del x,y\right>$ of the differential $\del$ on $CF_*(\phi)$ is then given by the mod 2 count of the index one components of the moduli space of holomorphic cylinders in $\RR\times M_\phi$, with respect to a generic $\RR$-invariant almost complex structure, which limit to the sections corresponding to $x$ and $y$ at $\pm\infty$, assuming these counts are finite. See \S\ref{floer} for details.

The standard condition, which in this form is due to Seidel \cite{s2}, used to ensure that these counts are finite is \emph{monotonicity}. Let $\omega_\phi$ be the two-form on $M_\phi$ induced by $\omega$ on $\Sigma\times\RR$ and let $c_\phi$ be the first Chern class of the vertical tangent bundle of $M_\phi \rightarrow S^1$.

\begin{defn}
\label{monointro}
A map $\phi: \Sigma \rightarrow \Sigma$ is \emph{monotone} if $[\omega_\phi] = k\,c_\phi$ as elements of $H^2(M_\phi; \RR)$ for some $k \in \RR$.
\end{defn}

Because $c_\phi$ controls the index, or expected dimension, of moduli spaces of holomorphic curves under change of homology class and $\omega_\phi$ controls their energy under change of homology class, this condition ensures that the energy is constant on the index one components of the moduli space, which implies compactness. See \S\ref{floerreview} and \S\ref{indexsection} for details.

The mapping class group $\Gamma_\Sigma$ is given by $\pi_0(\textrm{Diff}^+(\Sigma))$, the connected components of the group of oriented diffeomorphisms (or, equivalently, homeomorphisms) of $\Sigma$. For $h\in \Gamma_\Sigma$, let $\Symp_h(\Sigma)$ denote the space of symplectomorphisms in the mapping class $h$ and $\Symp_h^m(\Sigma)$ denote the space of monotone symplectomorphisms in the mapping class $h$. Seidel \cite{s2} has shown that the inclusions $$\Symp_h^m(\Sigma) \hookrightarrow \Symp_h(\Sigma) \hookrightarrow \Diff^+_h(\Sigma)$$ are homotopy equivalences; in particular, $\Symp_h^m(\Sigma)$ is path connected. Furthermore, Seidel showed that $HF_*(\phi)$ is invariant as $\phi$ is deformed through monotone symplectomorphisms. These imply that we have a symplectic Floer homology invariant $HF_*(h)$ canonically assigned to each mapping class $h$ given by $HF_*(\phi)$ for any monotone symplectomorphism.\footnote{Furthermore, when $\Sigma$ has negative Euler characteristic, these three groups, and in particular $\Symp_h^m(\Sigma)$, are contractible, so this assignment really is canonical.}

In this paper, we give an algorithm to compute $HF_*(h)$ in all previously unknown cases on surfaces of negative Euler characteristic.\footnote{The genus one case follows from Pozniak's thesis \cite[\S3.5.1-2]{po}.}
In addition, we define and give an algorithm to compute $HF_*(h)$ for all mapping classes on surfaces with boundary (see Remark \ref{dehnbdry} for our conventions for mapping classes on surfaces with boundary). We furthermore compute the $H_*(\Sigma)$-module structure. In a forthcoming paper \cite{c}, we extend the results of this paper to an algorithm to compute $HF_*(\phi)$ over a Novikov ring for any $\phi \in \Symp_h(\Sigma)$ and use this to give a sharp lower bound on the number of fixed points of an area-preserving map with nondegenerate fixed points in a given mapping class, generalizing the Poincar\'e-Birkhoff fixed point theorem.

Thurston's classification of surface diffeomorphisms \cite{t}, \cite{flp} states that every mapping class of $\Sigma$ is precisely one of
\bi
\ii Periodic (finite order): For some representative $\phi$, we have $\phi^\ell = \id$ for some $\ell \in \ZZ_{>0}$
\ii Pseudo-Anosov (see \S\ref{symprep})
\ii Reducible: Some representative $\phi$ fixes setwise a collection of curves $C$ none of which are nulhomotopic or boundary parallel (and the mapping class is not periodic)
\ei
In the reducible case, cutting along a maximal collection of curves $C$ gives a map on each component of $\Sigma\setminus C$ (given by the smallest power of $\phi$ which maps that component to itself) which is periodic or pseudo-Anosov.

\begin{remark}
\label{dehnbdry}
When $\Sigma$ has boundary, we instead consider the connected components of $\Diff^+_\del(\Sigma)$, the group of orientation-preserving diffeomorphisms of $\Sigma$ with no fixed points on the boundary. We then additionally consider maps with full Dehn twists at the boundary to be reducible, and periodic and pseudo-Anosov maps with boundary come with the additional data of which direction each setwise fixed boundary component rotates. See \S\ref{red} for more details.
\end{remark}

We briefly summarize previous work and which cases of $HF_*(h)$ have been computed. Dostoglou and Salamon \cite{ds} extended Floer's work \cite{ff} on symplectomorphisms Hamiltonian isotopic to the identity to other symplectomorphisms. Seidel \cite{s1} was the first to consider symplectic Floer homology specifically for surface symplectomorphisms, calculating it for arbitrary compositions of Dehn twists along a disjoint collection of curves. These are examples of reducible mapping classes, where the reducing collection of curves is the stated collection.
Gautschi \cite{g} calculated the symplectic Floer homology for all periodic mapping classes as well as reducible mapping classes in which the map on each component is periodic.

Eftekhary \cite{e} generalized Seidel's work on Dehn twists in a different direction. Given two disjoint forests $C_+$ and $C_-$ of embedded, homologically essential curves,\footnote{That is, any two curves in $C_+$ intersect at most once, and there are no cycles in the graph whose vertices are curves in $C_+$ and whose edges are intersections, and similarly for $C_-$. In addition, no curve in $C_+$ intersects any curve in $C_-$.} he showed that, for $\phi$ any composition of single positive Dehn twists along each curve in $C_+$ and single negative Dehn twists along each curve in $C_-$, we have $HF_*([\phi]) = H_{*(\mathrm{mod}\,2)}(\Sigma\setminus C_-, C_+)$. He also showed that some of these maps are pseudo-Anosov.

We determine the symplectic Floer homology for any pseudo-Anosov mapping class and for any reducible mapping class. In particular, we allow pseudo-Anosov components as well as periodic components.

For pseudo-Anosov mapping classes, we deform the canonical singular representative $\phi_{sing}$ of the mapping class (see \S\ref{symprep}) near its singularities to a smooth, symplectic representative with nondegenerate fixed points $\phi_{sm}$. For reducible mapping classes, we give a similar perturbation of the boundary behavior of the canonical singular representative of pseudo-Anosov maps in \S\ref{red} and combine this with standard representatives of maps on periodic components and of Dehn twists along the reducing curves to give standard-form-maps $\phi_{sm}$ for reducible mapping classes.

We give a notion of \emph{weak monotonicity} satisfied by $\phi_{sm}$ such that $HF_*(\phi)$ is well-defined for and invariant among weakly monotone maps. Monotonicity implies weak monotonicity, and so $HF_*(h) = HF_*(\phi)$ for any weakly monotone $\phi$. Furthermore we show that the $H_*(\Sigma)$-module structure of $HF_*(\phi)$ is invariant among weakly monotone maps.

\begin{defn}
\label{wkmonointro}
A map $\phi: \Sigma \rightarrow \Sigma$ is \emph{weakly monotone} if $[\omega_\phi]$ vanishes on $\ker(c_\phi|_{T(M_\phi)})$, where $T(M_\phi) \subset H_2(M_\phi;\RR)$ is generated by tori $T$ such that $\pi|_T: T \rightarrow S^1$ is a fibration with fiber $S^1$, where the map $\pi: M_\phi \rightarrow S^1$ is the projection.
\end{defn}

This brings us to our main results.

\begin{thm}
\label{pamain}For $h$ a pseudo-Anosov mapping class, we have:
\be
\ii
Every symplectomorphism $\phi:(\Sigma,\omega)\rightarrow(\Sigma,\omega)$ in the mapping class $h$ is weakly monotone. Hence if $\phi$ has nondegenerate fixed points, then $HF_*(\phi)$ is well-defined and $HF_*(\phi) = HF_*(h) = HF_*(\phi_{sm})$.
\ii Moreover, $HF_*(\phi_{sm}) = CF_*(\phi_{sm})$.
\ii Only the fundamental class of $H_*(\Sigma)$ (the identity), acts nontrivially.
\ee
\end{thm}

\begin{remark}
Furthermore, we show in \S\ref{tt} that $CF_*(\phi_{sm})$ is combinatorially computable from an invariant train track representation. In particular, see Theorem \ref{pacomb}. An invariant train track for a mapping class can be found algorithmically from, for example, a representation of the mapping class as a composition of Dehn twists by work of Bestvina and Handel \cite{bh}.
\end{remark}

\begin{thm}
\label{rmain}
For $h$ a reducible mapping class, we have:
\be
\ii Any standard form map is weakly monotone.
\ii $HF_*(h)$ splits into summands for each component\footnote{But note that the summand for a given component may be affected by the direction it is required to rotate near its boundary.} of $\Sigma\setminus C$ (and annular region; see \S \ref{red}), with the periodic and annular components as in Gautschi's paper \cite{g} and the pseudo-Anosov components as in Theorem \ref{pamain}, except when a pseudo-Anosov component with an unrotated boundary component (see \S\ref{red}) abuts a component on which the map is the identity, in which case there may be explicitly given additional flow lines between these two components. See Theorem \ref{redthm} for details and a complete statement and computation.
\ii $H_*(\Sigma)$ acts by intersection products on components on which the map is the identity and which do not meet pseudo-Anosov components, and only the fundamental class (the identity) acts nontrivially on other components.
\ee
\end{thm}

\begin{cor}
\label{growth}
The total rank of $HF_*(h^n)$
$$\left\{\begin{array}{ll}
\textrm{is bounded} & \textrm{if}\ \phi\ \textrm{is periodic} \\
\textrm{grows linearly} & \textrm{if}\ \phi\ \textrm{is reducible with all components periodic} \\
\textrm{grows exponentially} & \textrm{otherwise}
\end{array}\right.$$
\end{cor}

In \S\ref{exs} we explicitly work out $HF_*(h)$ for $h$ a large collection of pseudo-Anosov mapping classes based on a construction of Penner \cite{p} using Theorem \ref{pamain} and the train track computation in Theorem \ref{pacomb}:

\begin{thm}
\label{exmain}
Let $\{C_i\}$ and $\{D_j\}$ be collections of noncontractible simple closed curves on $\Sigma$ with the $C_i$ mutually disjoint, the $D_j$ mutually disjoint, and $C_i$ transverse to $D_j$. Suppose that every region of the complement $\Sigma-\coprod_i C_i - \coprod_j D_j$ is a disk with more than two edges.

Let $\phi = \prod_j \tau_{D_j}^{-2} \prod_i \tau_{C_i}^2$ and let $h$ be the mapping class of $\phi$. 
Then $HF_*(\phi) = HF_*(h)$ has total rank 
$$4\sum_{i,j} \left|C_i \cap D_j\right|^2 + 2g - 2.$$ 
\end{thm}

\subsection{Relationships with other Floer theories}

Symplectic Floer homology for surface symplectomorphisms is the $d = 1$ part of Periodic Floer Homology\footnote{Beware, however, that for some of the $\mathrm{spin}^c$-structures, the monotonicity condition for PFH differs from that for symplectic Floer homology.} (PFH) of $\phi$ \cite[\S2]{hs} \cite{h}, a Floer homology theory 
whose chain complex is generated by certain multisets of periodic orbits of $\phi$ and whose differential counts certain embedded pseudoholomorphic curves in $\RR\times M_\phi$. This theory was created to be a candidate for a $3$-dimensional version of Taubes' ``SW = Gr'' result \cite{ta1} and is conjectured \cite{hs} to be isomorphic to the Seiberg-Witten Floer homology (SWFH) \cite{km} of $M_\phi$ in certain $\textrm{spin}^c$-structures when $d < g-1$. This in turn is conjectured to agree with Ozsv\'ath and Szab\'o's Heegaard Floer homology invariant $HF^+$ \cite{os}.

PFH is also related to Embedded Contact Homology (ECH) \cite{hs2}, which is an analogous construction with $\RR \times Y$ for $Y$ a contact $3$-manifold
and is also conjectured to agree with SWFH and $HF^+$. Taubes has recently announced a result relating ECH to SWFH \cite{ta2, ta3} and yet more recently Lee and Taubes have announced a result relating PFH to SWFH. This would imply that the results of this paper calculate the SWFH of $M_\phi$ in certain $\mathrm{spin}^c$-structures.

PFH has been computed for a single Dehn twist and certain other cases \cite{hs}, but no pseudo-Anosov cases are known for $g > 1$. ECH has been computed for $T^3$ \cite{hs2} and for $M_\phi$ for $\phi$ an Anosov map on the torus (i.e. a map given by a hyperbolic element of $SL_2 \ZZ$) by Lebow \cite{leb}. Jabuka and Mark \cite{jm} have computed $HF^+(M_\phi)$ for a subset of the mapping classes calculated by Seidel \cite{s1} and Eftekhary \cite{e} for symplectic Floer homology, and the results agreed.

\subsection{Relationship with vanishing cycles and mutation}
Seidel has a conjecture for the symplectic Floer homology of a general composition of positive Dehn twists on a surface with boundary in terms of Hochschild Homology coming from his work on vanishing cycles and mutation \cite{s3}. 
Due to the complicated relationship between Thurston's classification of surface diffeomorphisms, which our results are based on, and compositions of positive Dehn twists, comparing Seidel's conjecture with the results of this paper is non-trivial. Since the first draft of this paper, Perutz \cite{per2} has announced a proof, in progress, of Seidel's conjecture. It would be interesting to find applications of the equality of these two very different computations. It may also be interesting to understand this equality directly, perhaps by understanding how the Floer complex changes as the Bestvina-Handel algorithm is implemented.

\subsection{Organization of the paper}

In \S\ref{floer} we review symplectic Floer homology and give a detailed discussion of when it is well-defined. We introduce Nielsen classes and the direct sum decomposition of $(CF_*,\del)$ resulting therefrom. We use a simplified version of Yi-Jen Lee's bifurcation analysis to show that, in our setting, $HF_*(\phi)$ is invariant under deformations through weakly monotone maps. Finally we show invariance of the $H_*(\Sigma)$-module structure using an algebraic interpretation of the $H_1(\Sigma)$-action in terms of the twisted Floer homology complex.

In \S\ref{wmpa} we show 
that every map in a pseudo-Anosov mapping class is weakly monotone.

In \S\ref{symprep}, we describe the canonical singular representative of a pseudo-Anosov mapping class $\phi_{sing}$ in the closed case and give a Hamiltonian perturbation supported near its singularities which results in a smooth symplectomorphism $\phi_{sm}$. In addition, we analyze the fixed points of $\phi_{sm}$ near the perturbed singularities and note that they are all of the same index.

In \S\ref{paproof}, we use a result of Birman and Kidwell \cite{bk} and the conclusions of \S\ref{symprep} to show that all differentials vanish on $CF_*(\phi_{sm})$ and deduce the $H_*(\Sigma)$-module structure.

In \S\ref{wmred}, we study the topology of the space of weakly monotone maps in arbitrary mapping classes in order to prove a stronger invariance result for $HF_*$.


In \S\ref{red}, we describe the standard form for reducible maps and discuss the behavior of pseudo-Anosov maps near the boundary. We give a Hamiltonian perturbation of $\phi_{sing}$ to $\phi_{sm}$ at the boundary. We then give in \S\ref{nielsenred} a classification of the Nielsen classes of fixed points of standard form reducible maps after Jiang and Guo \cite{jg}, extending their work to the symplectic case, where one must also consider fixed regions between parallel Dehn twists (which can be isotoped away if one is not required to preserve area). We then use this to show in \S\ref{redwkmono} that standard form reducible maps $\phi_{sm}$ are weakly monotone in every Nielsen class, which requires a lengthy technical argument. This allows us to show that $HF_*(\phi_{sm}) = HF_*(h)$.

In \S\ref{flred}, we use a neck stretching argument of Gautschi \cite{g} in Lemma \ref{cfsplits} to show that $(CF_*(\phi_{sm}),\del)$ splits into direct summands for each component (including regions between parallel Dehn twists) with one exception. This exception occurs when a pseudo-Anosov component abuts a component on which the map is the identity. We adapt Gautschi's neck stretching argument to this case and combine this with an argument of Salamon and Zehnder \cite{sz} in Lemma \ref{neck} to compute the differential by relating holomorphic cylinders from fixed points in one of these regions to the other to Morse flow-lines passing through the reducing curve between the regions. This allows us to give a computation in Theorem \ref{redthm} of $HF_*(\phi_{sm})$, together with its $H_*(\Sigma)$-module structure, using our result for pseudo-Anosov components and Gautschi's \cite{g} work on periodic components.

In \S\ref{tt}, we review the theory of train tracks and give a combinatorial formula computing $HF_*(\phi_{sm})$ for pseudo-Anosov mapping classes (and thus for reducible mapping classes as well, by the direct sum decomposition given in \S\ref{flred}) from the action of $\phi_{sm}$ on an invariant train track. Note that Bestvina and Handel \cite{bh} give an algorithm to find an invariant train track from e.g. the representation of a mapping class as a product of Dehn twists.

In \S\ref{exs}, we apply the results of \S\ref{tt} to a large collection of examples of pseudo-Anosov maps due to Penner and obtain an explicit formula for the rank of $HF_*(h)$ in these cases.

\subsection{Directions for further research}

\be
\ii \emph{Cylindrical contact homology}: Symplectic field theory and contact homology invariants \cite{egh} exist for mapping tori due to the existence of a Hamiltonian structure (see \cite{behwz, cm, fa}).
In this setting, the ``cylindrical mapping torus contact homology" splits as a direct sum $$\bigoplus_{n=1}^\infty H_*\left(C_*(\phi,n)\right),$$ where $C_*(\phi,n)$ is generated by good periodic orbits of degree $n$. This chain complex is similar to $CF_*(\phi^n)$, but the latter is generated instead by periodic \emph{points}, with multiplicity, of degree $n$. This, and its more complicated ``full'' version, are a sort of intermediate step between symplectic Floer homology and (cylindrical) contact homology of open books. Mei-Lin Yau \cite{y} has given a correspondence between (not necessarily cylindrical) holomorphic curves in mapping tori and in their $0$-surgeries, which are contact $3$-manifolds with a natural open book decomposition. Since the first draft of this paper came out, Colin and Honda \cite{ch} have computed certain examples of cylindrical contact homology for open books.

\ii \emph{Product structures} We have an understanding of the $H_*(\Sigma)$-module structure, but lack an understanding of products $*: HF(\phi)\otimes HF(\psi)\rightarrow HF(\phi\circ\psi)$ unless one of the two maps is the identity. An important special case is an understanding of rings of the form $$\bigoplus_{n=0}^\infty HF_*(\phi^n) \ \mathrm{or}\ \bigoplus_{n=-\infty}^\infty HF_*(\phi^n).$$

\ii \emph{Grading and $\mathrm{spin}^c$-structures} We use a $\ZZ/2$-grading, but a relative $\ZZ/(2\mathrm{div}(c_\phi)) = \ZZ/(4g-4)$ grading is available for each Nielsen class\footnote{In our calculations, we have more than one generator in a given Nielsen class only in situations in which we have a correspondence between Floer flow-lines and Morse flow-lines. Thus the relative grading for Nielsen classes is simply the Morse relative grading.} or more generally for each homological Nielsen class. That is, among generators whose orbits agree in $H_1(M_\phi)$. A good understanding of the splitting not into Nielsen classes but into homological Nielsen classes, together with relative grading data, is lacking.

\ee

\subsection*{Acknowledgments} I'd like to thank my advisor Michael Hutchings for his invaluable guidance and support and David Farris, Eli Lebow, and Lee Mosher for helpful discussions. I'd also like to thank an anonymous referee for very thorough and insightful comments.

\section{Floer homology, weak monotonicity, Nielsen classes, bifurcation analysis, and module structure}
\label{floer}

\subsection{Review of symplectic Floer homology}
\label{floerreview}

Let $\Sigma$ be a compact, connected, oriented surface of negative Euler characteristic. Let $\omega$ be a symplectic form (i.e. an area form) on $\Sigma$. Let $\phi \in \Symp(\Sigma,\omega)$ be a symplectomorphism (i.e. an area-preserving diffeomorphism) of $\Sigma$ and consider $$M_\phi = \frac{\RR\times\Sigma}{(t+1,x) \sim (t,\phi(x))},$$ the mapping torus of $\phi$. Note that this is a $\Sigma$-bundle over $S^1$ and we have the projection $\pi: M_\phi \rightarrow \RR/\ZZ = S^1$.

Suppose $\phi$ has nondegenerate fixed points; that is, at a fixed point, $d\phi$ does not have $1$ as an eigenvalue (in particular, fixed points are isolated). Under certain conditions (see below), we can define the symplectic Floer homology $HF_*(\phi)$ associated to $\phi$. This is the homology of a $\ZZ/2$-graded\footnote{There is a lift to a relative $\ZZ/\mathrm{div}(2c_1(E)) = \ZZ/(4g-4)$ grading, with a corresponding absolute grading by homotopy classes of $2$-plane fields; see \cite{per1} or \cite{h3}.} chain complex $CF_*(\phi)$ (which depends on the choice of almost complex structure; see below) which is freely generated as a vector space over $\ZZ/2$ by the fixed points of $\phi$. The grading $\epsilon(x)$ of a fixed point $x$ is given by $(-1)^{\epsilon(x)} = sign\left(\det(1-d\phi)\right)$. The homology, but not necessarily the chain complex itself, is invariant under certain deformations of $\phi$ (and independent of the choice of almost complex structure). We give a summary of how to define $HF_*(\phi)$ and point out any relevant pitfalls. See also \cite{ds, s2}.

Let $\Gamma(M_\phi)$ denote the space of smooth sections of $\pi: M_\phi \rightarrow S^1$, also known as the twisted loop space,\footnote{We will abuse notation and switch between sections of $M_\phi$ and paths in $\Sigma$ at will.} consisting of smooth paths $\gamma: \RR \rightarrow \Sigma$ such that $\gamma(t) = \phi(\gamma(t+1))$. Note that a fixed point $x \in \Sigma$ of $\phi$ can be interpreted as a constant path $\gamma(t) = x$ in $\Gamma(M_\phi)$. The tangent space of $\Gamma(M_\phi)$ at a path $\gamma$ consists of sections $\xi$ of $\gamma^*T\Sigma$ such that $\xi(t) = d\phi(\xi(t+1))$. We have an action one-form $\alpha_\phi$ on $\Gamma(M_\phi)$ defined by $$\alpha_\phi(\gamma) \xi = \int_0^1 \omega(\frac{d\gamma}{dt}, \xi)\ dt.$$ Note that the zeroes of this one-form are precisely the constant paths $\gamma$, that is, the paths corresponding to fixed points of $\phi$.

The differential on $CF_*(\phi)$ is defined by a certain count of ``gradient flow lines'' of $\alpha_\phi$. To get a metric (which doesn't show up explicitly in the equations), we require a complex structure on $\Sigma$. Let $\mathcal{J}$ be the space of smooth $1$-periodic families $J = \{J_t\}$ (that is, $\phi_*J_{t+1} \phi_*^{-1}= J_t$) of $\omega$-tame\footnote{That is, $\omega(v,Jv) > 0$ for $v \neq 0$.} almost complex structures on $\Sigma$. Then we can define $\MM(x,y,J)$ to be the set of maps $u: \RR^2 \rightarrow \Sigma$
such that 
\beqan
&&u(s,t) = \phi(u(s,t+1)) \label{d1}\\
&&\lim_{s\rightarrow +\infty} u(s,\cdot) = x\ \ \textrm{and}\ \ \lim_{s\rightarrow -\infty} u(s,\cdot) = y \\
&&\del_s u + J_t(u) \del_t u = 0 \label{d3}
\eeqan
We call such a $u$ a \emph{flow line} from $x$ to $y$. We may think of $u$ as a section $u: \RR\times S^1 \rightarrow \RR\times M_\phi$ of $\id\times \pi: \RR\times M_\phi \rightarrow \RR\times S^1$. We may also think of $u$ as a path in $\Gamma(M_\phi)$ between $\gamma_x$ and $\gamma_y$, the constant paths associated to the fixed points $x$ and $y$. We will use these three interpretations interchangeably.

At each flow line $u$, we have a Fredholm linearized $\dbar$ operator $D_u$.
We let $\MM_k(x,y,J)$ be the subset consisting of those maps $u$ for which the index of $D_u$ is $k$. For generic $J \in \mathcal{J}$, $D_u$ is onto and so by the implicit function theorem, all of these are smooth manifolds of dimension $k$ (see \cite{fhs}). Additionally, these moduli spaces come with a natural action of $\RR$ by translation in the first coordinate, and so we consider their quotients $\MM_k(x,y,J)/\RR$.

We would like to define the differential\footnote{We use symplectic Floer \emph{homology}, which agrees with the conventions in Periodic Floer Homology \cite{h} \cite{hs}.}
 by \beqan\label{diff}\del x = \sum_y c_{x,y}\ y\eeqan where $c_{x,y}$ is a count of the number of points (modulo 2) in $\MM_1(x,y,J)/\RR$, and we would also like to have $\del^2 = 0$ so that we can take the homology $HF_*(\phi)$. In order to do this, we need compactness and gluing results (note that bubbling is not an issue as $\pi_2(\Sigma) = 0$), which require upper bounds on the energy 
\beqan\label{e1} E(u) = \frac{1}{2}\int_{\RR\times [0,1]} \left(\left|\del_s u\right|^2 + \left|\del_t u\right|^2\right) ds\,dt\eeqan in each $\MM_k(x,y,J)$ (at least for $k = 1$ or $2$). For $u \in \MM_k(x,y,J)$, we have \beqan\label{e2} E(u) = \int_{\RR\times [0,1]} u^*\omega. \eeqan

These are locally constant. To ensure they are bounded, we impose:\footnote{We could make do without this if we were willing to work over a Novikov ring. This is the approach we take in \cite{c}.}

\begin{cond}
\label{econd}
$E(u)$ is constant on $\MM_k(x,y,J)$.\end{cond}

Given this, the symplectic Floer homology groups $HF_*(\phi)$ are well defined. Additionally, they are independent of the choice of (generic) $J$ (which follows, for example, by Theorem \ref{invariance}).
We have yet to consider invariance under (certain) deformations of the map $\phi$. More on that in \S\ref{bif}.

\subsection{Index and weak monotonicity}
\label{indexsection}

In order to verify Condition \ref{econd}, we need to know more about the index $k = \ind(D_u)$. We follow \cite[\S2.1-2.3,5.1]{h}. We first define the relative first Chern class and the Conley-Zehnder index and then use these to state a formula for $\ind(D_u)$.

Let $E$ denote the vertical tangent bundle of $M_\phi \rightarrow S^1$. This bundle carries a symplectic structure induced by $\omega$. If $\gamma_x \in \Gamma(M_\phi)$ corresponds to a fixed point $x$ (i.e. $\gamma_x(t) = x$), let $\TT(\gamma_x)$ denote the set of homotopy classes of symplectic trivializations of $\gamma_x^* E$. This is an affine space over $\ZZ$.

Suppose we have a flow line $u: \RR\times S^1 \rightarrow \RR\times M_\phi$ which limits to $x$ at $-\infty$ and $y$ at $+\infty$. If we choose trivializations $\tau_x \in \TT(\gamma_x)$ and $\tau_y \in \TT(\gamma_y)$, we get a \emph{relative first Chern class} $c_1(u^*E,\tau_x,\tau_y) \in \ZZ$ by counting the zeroes of a section of $u^*E$ which is constant (and non-zero) with respect to $\tau_x$ and $\tau_y$ over the two boundary circles. If $v$ is another flow line with the same limits, let $[u-v]$ be the homology class in $H_2(\RR\times M_\phi) = H_2(M_\phi)$ corresponding to the obvious (up to reparametrization) map $S^1\times S^1 \rightarrow \RR\times M_\phi$. Then we have
\beqan \label{relcc} c_1(u^*E,\tau_x,\tau_y) - c_1(v^*E,\tau_x,\tau_y) = \left<c_1(E), [u-v]\right>.\eeqan

We also have an integer called the \emph{Conley-Zehnder index} $\mu(\gamma_x,\tau_x)$ associated to a fixed point $x$ and a trivialization $\tau_x \in \TT(\gamma_x)$. We have a natural connection on the bundle $M_\phi \rightarrow S^1$ whose monodromy is given by $\phi$. Linearized parallel transport induces a connection on $\gamma_x^*E$. Thus if we traverse $\gamma_x$ once, then using $\tau_x$ this induces a path of symplectic matrices in $\textrm{Sp}(2,\RR)$ from $1$ to $d\phi_x$. We let $\mu(\gamma_x,\tau_x)$ be the Maslov index of this path of symplectic matrices (see e.g. \cite{sz}). 
Furthermore, $\mu(\gamma_x,\tau_x) \equiv \epsilon_x$ (modulo 2).

\begin{thm}[{Index, \cite{rs}}]
\label{index}
For $u$ a flow line from $x$ to $y$, $$\ind(D_u) = 2c_1(u^*E,\tau_x,\tau_y) + \mu(\gamma_x,\tau_x) - \mu(\gamma_y,\tau_y).$$ In particular, this is independent of the choice of trivializations.\footnote{We use the notation of Hutchings \cite[Theorem 5.1]{h}, who states this theorem for general Riemann surfaces, where it is due to Schwarz \cite[Theorem 3.3.11]{sc}.}\qed
\end{thm}

\begin{cor}[Change of homology class]
\label{chindex}
If $u$ and $v$ are flow lines with the same limits, then $$\ind(D_u) - \ind(D_v) = 2\left<c_1(E),[u-v]\right>.$$\qed
\end{cor}

Following Seidel \cite{s2}, we now turn this into a way to verify Condition \ref{econd}. If $u$ and $v$ are two flow lines from $x$ to $y$, then by following $u$ and then the reverse of $v$, we get a loop in $\Gamma(M_\phi)$ based at $\gamma_x$, i.e. an element $\psi_{uv} \in \pi_1(\Gamma(M_\phi), \gamma_x)$. We may think of an element $\psi \in \pi_1(\Gamma(M_\phi),\gamma_x)$ as a map $S^1\times S^1 \rightarrow S^1 \times M_\phi$ which is a section of $\id \times \pi: S^1 \times M_\phi \rightarrow S^1 \times S^1$. On $S^1 \times M_\phi$ as well as on $M_\phi$, we have two elements of $H^2(M_\phi; \RR)$ (we will abuse notation and consider these on $S^1\times M_\phi$ or $M_\phi$ interchangeably). First we have $\omega_\phi$ induced by $\omega$ on $\Sigma$, and second we have $c_\phi = c_1(E)$. Both of these may be evaluated on elements of $\pi_1(\Gamma(M_\phi),\gamma_x)$. We denote these actions by $\omega_\phi(\psi) = \int_{S^1\times S^1} \psi^*\omega_\phi$ and $c_\phi(\psi) = c_1(\psi^*E)[S^1\times S^1]$.

\begin{cor}[Monotonicity for $x$, {cf. \cite[Lemma 9]{s2}}]
\label{echeck}
Suppose that for all $\psi\in \pi_1(\Gamma(M_\phi),\gamma_x)$ such that $c_\phi(\psi) = 0$, we also have $\omega_\phi(\psi) = 0$.
Then $E(u)$ is constant on $\MM_k(x,y,J)$ for all $k$, $y$, and $J$. That is, Condition \ref{econd} is satisfied for flow lines starting at $x$. In this case, we say that $\phi$ satisfies \emph{monotonicity for $x$}.
\end{cor}

\proof
Suppose we have $u$ and $v$ in $\MM_k(x,y,J)$. Then we have $\ind(D_u) - \ind(D_v) = k - k = 0$, so  by Corollary \ref{chindex}, $c_\phi(\psi_{uv}) = \left<c_1(E),[u-v]\right> = 0$. Thus, using $\omega_\phi(\psi_{uv}) = 0$, we have:
 $$0 = \omega_\phi(\psi_{uv}) = \int_{S^1\times S^1} \psi^*\omega_\phi = \int_{\RR\times[0,1]} u^*\omega - \int_{\RR\times[0,1]} v^*\omega = E(u) - E(v).$$
\qed

Thus we've reduced Condition \ref{econd} to:

\begin{cond} [Weak Monotonicity]
\label{wcond}
$\omega_\phi$ vanishes on the kernel of $c_\phi: H_1(\Gamma(M_\phi)) \rightarrow \RR$
\end{cond}

\subsection{Nielsen classes}
\label{nclass}

Given a map $\phi$, there is a topological separation of fixed points due to Nielsen. The chain complex $CF_*(\phi)$ splits into direct summands corresponding to each Nielsen class (defined below). That is, if two fixed points $x$ and $y$ are in different Nielsen classes, then the matrix coefficient $c_{x,y} = \left<\del x,y\right>$ is zero. This idea has been used e.g. by Pozniak \cite[\S 3.3.1, 3.5.1]{po} and Seidel \cite{s1}, and has also been noticed by Nielsen theorists, e.g. Fel'shtyn \cite{fe}.

\begin{defn} Two fixed points $x$ and $y$ of $\phi$ are in the same \emph{Nielsen class}, or \emph{Nielsen equivalent}, if there exists a path $\gamma: [0,1] \rightarrow \Sigma$ with $\gamma(0) = x$, $\gamma(1) = y$ such that $\phi(\gamma)$ is homotopic to $\gamma$ (rel. endpoints). Note that this is an equivalence relation.\end{defn}

There can be a flow line from $x$ to $y$ only if $x$ and $y$ are in the same Nielsen class: If $u(s,t)$ is a flow line from $x$ to $y$, let the path $\gamma_t : [0,1] \rightarrow \Sigma$ be a reparametrization of $u(.,t)$ (for each $t$). Then the curve $\gamma_1$ exhibits the Nielsen equivalence between $x$ and $y$: its limits are $x$ and $y$, and it is homotopic (rel. endpoints) through $\gamma_t$ to $\gamma_0 = \phi(\gamma_1)$.

Interpreting fixed points as sections of $M_\phi \rightarrow S^1$, that is, as elements of $\Gamma(M_\phi)$, we note that two fixed points are Nielsen equivalent if and only if their sections are homotopic. That is, $\pi_0(\Gamma(M_\phi))$ describes possible Nielsen classes.

For $\eta \in \pi_0(\Gamma(M_\phi))$, let $CF_*(\phi,\eta)$ be the chain complex generated by fixed points in Nielsen class corresponding to $\eta$. (Note that this may be the zero chain complex.) We have $CF_*(\phi) = \bigoplus_{\eta\in \pi_0(\Gamma(M_\phi))} CF_*(\phi,\eta)$ as chain complexes, and thus $$HF_*(\phi) = \bigoplus_{\eta\in \pi_0(\Gamma(M_\phi))} HF_*(\phi,\eta).$$

Let $\Gamma(M_\phi)_\eta$ denote the $\eta$-component of $\Gamma(M_\phi)$. Then $HF_*(\phi,\eta)$ is well-defined if we have:

\begin{cond} [$\eta$-Weak Monotonicity]
\label{ncond}
$\omega_\phi$ vanishes on the kernel of $c_\phi: H_1(\Gamma(M_\phi)_\eta) \rightarrow \RR$
\end{cond}

\subsection{Bifurcation analysis of Yi-Jen Lee}
\label{bif}

The standard argument to show invariance under certain deformations of $\phi$ is the continuation argument as in \cite{ff}. It gives invariance under Hamiltonian deformations of $\phi$. Seidel \cite[\S 3]{s2} has a generalized continuation argument which shows invariance under deformations which preserve monotonicity. We have need of a more general version of monotonicity, and thus also require a more general invariance result.

Floer's original invariance argument \cite{fl} was not a continuation argument, but rather a bifurcation argument, in which one analyzes degeneracies that arise in generic one-parameter families and their effect on $HF_*$. He did not complete all of the analysis required, but with the advent of the continuation argument, this was no longer needed. Yi-Jen Lee \cite{l1,l2} has worked out a general bifurcation argument for what she calls Floer-type theories, including all of the necessary analysis. 
Michael Usher \cite{u} has a nice summary of the invariance result this gives (which Lee conjectured in an earlier paper \cite[Eqn 3.2]{l} but did not explicitly state as a Theorem in \cite{l1,l2}) and its algebraic aspects. We give a brief discussion of a simplified version sufficiently general for our needs.

\begin{thm}[Simplified version of {\cite[Theorem 3.6, due to Lee (\cite{l1},\cite{l2})]{u}}]
\label{invariance}
Suppose $(\Sigma, \omega)$ is a symplectic manifold with $\pi_2(\Sigma) = 0$. Let $\phi_r: \Sigma \rightarrow \Sigma$ be a smooth family of symplectomorphisms and $J_r = \{J_t\}_r$ a smooth family (of $1$-periodic families) of almost complex structures such that $(\phi_0,J_0)$ and $(\phi_1,J_1)$ are generic. Let $\eta \in \pi_0(\Gamma(M_{\phi_0})) \cong \pi_0(\Gamma(M_{\phi_r}))$ be a Nielsen class. Suppose that $\omega_{\phi_r}$ vanishes on $\ker(c_\phi): \pi_1(\Gamma(M_{\phi_r}), \eta) \rightarrow \RR$ (that is, $\phi_r$ is $\eta$-weakly monotone).
Then $(CF_*(\phi_0,\eta, J_0), \del_0)$ is chain homotopy equivalent to $(CF_*(\phi_1,\eta, J_1), \del_1)$.\qed
\end{thm}

The main idea of the proof is to study generic $1$-parameter families and what happens to $(CF_*, \del)$ at times $r_0$ where the data fail to be generic. Because we have assumed that $\omega_\phi$ vanishes on $\ker(c_\phi)$, we cannot have a flow line of index $0$ from a fixed point $x$ to itself (Lee calls these ``type II handleslides'') because, by Corollary \ref{chindex}, it would have zero energy.
The argument is then formally identical to bifurcation arguments in Morse theory; see for instance \cite{h2}
for a carefully done modern perspective in the more complicated circle-valued case (which includes type II handleslides). For Lee's account in the Floer theory setting, see \cite[\S4.4]{l1}, noting that we have no type II handleslides.

In generic $1$-parameter families, there are three types of degeneracies that can happen at a time $r_0$:

\bi
\ii \textbf{Cancellation}: A degenerate flow-line from $x$ to $y$ of index 1. 
\ii \textbf{Birth-death}: A birth or death of two fixed points $x$ and $y$ at a degenerate fixed point.
\ii \textbf{Handleslide}: A flow-line from $x$ to $y$ of index 0.
\ei

For a cancellation, at times $r_0 - \epsilon < r < r_0$ we have two flow-lines (of opposite sign were we working with $\ZZ$-coefficients) which cancel at time $r_0$ (or the time-reversal of this situation). This does not change the chain complex $(CF_*, \del)$.

For a birth-death, let us assume without loss of generality we have a death of a pair $x$ and $y$, with a single flow-line from $x$ to $y$ for $r$ shortly before $r_0$. Let $z$ and $w$ be two other fixed points. We compare $\left<\del_- z, w\right>$ and $\left<\del_+ z,w\right>$, the matrix coefficients at times just before $r_0$ and at times just after $r_0$, respectively. Because $\omega_\phi$ vanishes on $\ker(c_\phi)$, we may have either flow-lines of index $1$ from $z$ to $y$ or flow-lines of index $1$ from $y$ to $z$, but not both. Let us assume that the former is the case. Flow lines from $z$ to $y$ are glued to flow lines from $x$ to $w$. Thus $$\left<\del_+ z,w\right> = \left<\del_- z, w\right> \pm \left<\del_- z, y\right>\left<\del_- x, w\right>.$$ This does not change the chain homotopy type of the chain complex (we need not worry about the sign as we are working with $\ZZ/2$-coefficients, but the case with $\ZZ$-coefficients is fine as well); see e.g. \cite[\S 3.3]{u}.

For a handleslide, the generators of the chain complex do not change. We again compare $\del_-$ and $\del_+$. If we have a fixed point $z$ with a flow-line to $x$ for $r$ near $r_0$, then, by gluing this with our flow-line from $x$ to $y$, we get a flow-line from $z$ to $y$ on one side of $r_0$. Similarly, if we have a fixed point $w$ with a flow-line from $y$ for $r$ near $r_0$, then we get a flow-line from $x$ to $w$ on one side of $r_0$. Thus, letting $A: CF_* \rightarrow CF_*$ be given by $A(x) = x \pm y$ and $A(v) = v$ for all $v \neq x$, we have that $\del_+ = A^{-1}\circ \del_- \circ A$. Thus the isomorphism-type of $(CF_*,\del)$ does not change, and, in particular, the homology is unchanged.

The analysis required to make the above rigorous is the gluing described in the birth-death situation and in the handleslide situation. Lee does this in the second paper \cite{l2}, where she must also glue tori due to type II handleslides.

\subsection{Product structure and invariance}
\label{product}

There is a product $*: HF(\phi)\otimes HF(\psi) \rightarrow HF(\phi\circ\psi)$, whose matrix coefficients count holomorphic sections of a hamiltonian symplectic fibration with fiber $\Sigma$ over the pair of pants $\CC\PP^1 \setminus \{0,1,\infty\}$. The monodromies are given by $\phi$, $\psi$, and $\phi\circ\psi$.

In particular, this gives $HF(\phi)$ the structure of a module over the quantum homology $QH(\Sigma) = HF(id)$ via $*: HF(id)\otimes HF(\phi) \rightarrow HF(\phi)$. By Liu and Tian \cite{lt}, this intrinsic definition is equivalent to the following extrinsic definition: given $a \in HF(id)$ and x a fixed point of $\phi$, choose a cycle $C$ representing $a \in HF(id) \cong H_*(\Sigma)$. Then the coefficient $\left<a*x,y\right>$ for $y$ a fixed point of $\phi$ is given by a count holomorphic cylinders as in Section \ref{floerreview} with the extra condition that they pass through the cycle $C$ at $s = t = 0$.

The action of the fundamental class on $HF(\phi)$ is the identity. The action of the point class has been studied by Seidel \cite{s2}, who shows that it acts as zero unless the mapping class of $\phi$ is the identity mapping class. We calculate the action of $H_1(\Sigma)$ below. For this to have meaning, we have need of a proof of the invariance of the module structure under weakly monotone deformations, which we come to presently.

\begin{prop}
\label{productinvariant}
The action of $H_1(\Sigma)$ on $HF(\phi)$ is invariant under weakly monotone deformations of $\phi$.
\end{prop}

\begin{remark}
This implies that the full module structure is invariant by the discussion in the previous paragraph, or by the observation that the action of $H_1(\Sigma)$ algebraically determines the full module structure. We expect the full module structure to be invariant for weakly monotone deformations of symplectomorphisms in higher dimensions, using more details of Lee's bifurcation analysis \cite{l1,l2}.
\end{remark}

\proof
The key observation here is that the action of $H_1(\Sigma)$ has index one and thus, using the extrinsic definition, we may consider index one moduli spaces of cylinders, the same moduli spaces used in defining the differential, and then impose the constraint that they pass through $C$ at $s = t = 0$ (whereas for the differential we mod out by the $\RR$-action). We fix $C$ missing the fixed points of $\phi$. Then we have
\beqan\label{h1action}\left<[C]*x,y\right> = \sum_{p \in \MM_1(x,y,J)/\RR} [C]\cdot[p],\eeqan where $[p]$ is the relative homology class of the cylinder from $x$ to $y$, in the $H_2(M_\phi)$-torsor of $2$-chains with boundary $\gamma_y - \gamma_x$, and $[C]\cdot[p]$ is the intersection product, where by $[C]$ we really mean $i_*[C]$ for $i:\Sigma\rightarrow M_\phi$ the inclusion to the fiber over $0 \in S^1$.

The extrinsic definition is equivalent to the following algebraic definition, as in \cite[\S 12.1.3]{hs2}. Let $H = H_2(M_\phi)$. We can define the fully twisted symplectic Floer homology $(CF^{tw}(\phi),\del^{tw})$ as the free $\ZZ/2[H]$-module generated by fixed points of $\phi$, with differential given by $$\left<\del^{tw} x,y\right> = \sum_{p \in \MM_1(x,y,J)/\RR} ([p]-b_{xy}),$$ where $b_{xy}$ is a choice of an element in the $H$-torsor associated to $x$ and $y$. We require the choices to be compatible in the sense that they form a cochain with respect to the gluing maps and additionally require that $[C]\cdot b_{xy} = 0$ for all $x,y$. This is possible because we may assume that $C$ is nonseparating (else $[C]$ is zero). The chain homotopy type of $(CF^{tw}(\phi),\del^{tw})$ is invariant under weakly monotone deformations of $\phi$. The proof is the same as that of Theorem \ref{invariance}, with the additional comment that gluing is additive on relative homology classes, and so the same algebraic formulas given above for birth-deaths and handleslides carry over.\footnote{For the handleslide, we actually must use $A(x) = x \pm (h-b_{xy}) y$, where $h$ is the class of the index zero flow-line. For the birth-death, we simply choose $b_{xy}$ to be the class of the short flow-line between $x$ and $y$ and the same formula carries over.}

Via the intersection product, the chain $C$ gives a map $i_C: H \rightarrow \ZZ/2$. This induces a $\ZZ/2$-linear map $\zeta_C: \ZZ/2[H] \rightarrow \ZZ/2[H]$ defined on elements $h\in H$ by $\zeta_C(h) = i_C(h)h$. Note that this is a derivation of $\ZZ/2[H]$ over $\ZZ/2$. We now define a $\ZZ/2[H]$-linear map $\del^{tw}_C$ on $CF^{tw}(\phi)$ via $$\del^{tw}_C x = \sum_y \zeta_C\left(\left<\del^{tw} x, y\right>\right)y.$$

We have that $\del^{tw} \circ \del^{tw}_C = \del^{tw}_C \circ \del^{tw}$ (as noted in \cite[\S 12.1.3]{hs2}, this follows from $(\del^{tw})^2 = 0$ and $\zeta_C(xy) = \zeta_C(x)y + x\zeta_C(y)$). Setting all elements of $H$ equal to one, we have $\del_C$ on $CF(\phi)$ which commutes with $\del$. This gives a map $\del_C:HF(\phi)\rightarrow HF(\phi)$ which agrees with the map $[C]*-:HF(\phi)\rightarrow HF(\phi)$ by equation (\ref{h1action}). To complete the proof of the invariance of the latter under weakly monotone deformations, we have only to show that the former is invariant under chain homotopies of $CF^{tw}(\phi)$.

Consider a $\ZZ/2[H]$-linear chain map $\Phi: CF^{tw}(\phi) \rightarrow CF^{tw}(\phi')$ (we're only interested in the case where this is a quasi-isomorphism, but it does not matter). Note that via the inclusion $i: CF \rightarrow CF^{tw}$ and the map $j: CF^{tw} \rightarrow CF$ sending all elements of $H$ to one, we have a chain map $j\circ\Phi\circ i: CF(\phi) \rightarrow CF(\phi')$. We show $\Phi \del^{tw}_C$ and $\del^{tw}_C \Phi$ are chain homotopic, which implies $j\Phi i \del_C$ and $\del_C j\Phi i$ are chain homotopic. We compute:
\beqa
\left<\Phi \del^{tw}_C x_i, x_k\right> &=& \left<\Phi\left<\del^{tw}_C x_i, x_j\right> x_j, x_k\right> \\
&=& \left<\del^{tw}_C x_i, x_j\right>\left<\Phi x_j,x_k\right> \\
&=& \zeta_C\left(\left<\del^{tw} x_i,x_j\right>\right)\left<\Phi x_j,x_k\right> \\
&=& \zeta_C\left(\left<\del^{tw} x_i,x_j\right>\left<\Phi x_j,x_k\right>\right) - \left<\del^{tw} x_i, x_j\right> \zeta_C\left(\left<\Phi x_j,x_k\right>\right)
\eeqa

Similarly, \beqa \left<\del^{tw}_C \Phi x_i,x_k\right> &=& \zeta_C\left(\left<\Phi x_i, x_j\right>\left<\del^{tw}x_j,x_k\right>\right) - \zeta_C\left(\left<\Phi x_i,x_j\right>\right)\left<\del^{tw}x_j,x_k\right>
\eeqa

The first terms on the right hand sides are equal because $\del^{tw}\Phi = \Phi \del^{tw}$. The others are of the form $\del^{tw}K$ and $K\del^{tw}$, so indeed the two left hand sides are chain homotopic.

\qed

\section{Pseudo-Anosov maps}

\subsection{Weak monotonicity and invariance in pseudo-Anosov mapping classes}
\label{wmpa}

The aim of this section is to show:

\begin{prop}
\label{notmono}
For any (not necessarily monotone) symplectomorphism $\phi:(\Sigma,\omega)\rightarrow(\Sigma,\omega)$ in a pseudo-Anosov mapping class $h$, if $\phi$ has nondegenerate fixed points, then $HF_*(\phi)$ is well-defined and $HF_*(\phi) = HF_*(h)$.
\end{prop}

In order to prove this we show that every symplectomorphism $\phi$ in a pseudo-Anosov mapping class is weakly monotone (i.e. satisfies Condition \ref{wcond}). To do this, we show that $H_1(\Gamma(M_\phi)) = 0$.

\begin{lemma}
\label{pi1is0}
For $\phi$ pseudo-Anosov and $\gamma \in \Gamma(M_\phi)$, $\pi_1(\Gamma(M_\phi), \gamma) = 0$.
\end{lemma}

\proof
We consider $\gamma$ as a map $\RR\rightarrow\Sigma$ with $\gamma(t) = \phi(\gamma(t+1))$. An element of $\pi_1(\Gamma(M_\phi),\gamma)$ is of the form $\gamma_s(t)$ for $s\in S^1 = \RR/\ZZ$ with $\gamma_0(t) = \gamma(t)$.
We consider $\alpha_0(s) = \gamma_s(0)$ and $\alpha_1(s) = \gamma_s(1)$. These are closed curves on $\Sigma$ and $\phi(\alpha_1(s)) = \alpha_0(s)$. Furthermore, $\alpha_0(s)$ is homotopic to $\alpha_1(s)$ by the homotopy $\alpha_t(s) = \gamma_s(t)$. We now quote Lemma \ref{nohom} (below) to conclude that $\alpha_0(s)$ must be nulhomotopic.

To deal with this case, we consider the fibration 

$$\begin{CD}
\Omega\Sigma @>>> \Gamma(M_\phi) \\
@. @VVfV \\
@. \Sigma \\
\end{CD}.$$
Here $f$ takes a section to its value at $t = 0$. The fiber over $\gamma(0) \in \Sigma$ is homeomorphic to $\Omega\Sigma$, which we think of as loops in $\Sigma$ based at $\gamma(0)$, as follows: Let $\psi_t : \Sigma \rightarrow \Sigma$ be an isotopy from the identity to a map which takes $\gamma(0)$ to $\gamma(1)$. Then, given an element $x(t)$ of $\Omega\Sigma$ (i.e. a loop based at $\gamma(0)$), we consider $y(t) = \psi_t(x(t))$, an element of the fiber over $\gamma(0)$ of the map $f$. This gives the homeomorphism.

Thus we have a long exact sequence of homotopy groups, a piece of which is \beqa \pi_2(\Sigma,\gamma(0)) \rightarrow \pi_1(\Gamma(M_\phi), \gamma) \rightarrow \pi_1(\Sigma,\gamma(0)).\eeqa The image in $\pi_1(\Sigma,\gamma(0))$ of the element of $\pi_1(\Gamma(M_\phi),\gamma)$ represented by the homotopy $\gamma_s(t)$ is represented by $\alpha_0(s)$, which is nulhomotopic. Thus it comes from $\pi_2(\Sigma,\gamma(0)) $, which is trivial. Thus we have $\pi_1(\Gamma(M_\phi),\gamma) = 0$.
\qed

\begin{lemma}[{\cite[Lemma 5.6]{cb}}]
\label{nohom}
If $\phi:\Sigma\rightarrow\Sigma$ is pseudo-Anosov, and $C$ is a noncontractible curve on $\Sigma$, then $\phi^m(C) \simeq \phi^n(C)$ if and only if $m = n$. \qed
\end{lemma}

\textit{Proof of Proposition \ref{notmono}:}\ \ \ \ 
By Lemma \ref{pi1is0}, Condition \ref{wcond} is
trivially satisfied. 
Thus $HF_*(\phi)$ is well-defined. To show that $HF_*(\phi) = HF_*(h)$, we use Theorem \ref{invariance}.
We find an isotopy of symplectic maps $\phi_t$ such that $\phi_0 = \phi$ and $\phi_1$ is monotone. This exists because $\Symp^m \hookrightarrow \Symp$ is a homotopy equivalence \cite{s2}. By Lemma \ref{pi1is0}, each $\phi_t$ satisfies the conditions of Theorem \ref{invariance}.\qed

\subsection{Symplectic representatives of pseudo-Anosov mapping classes}
\label{symprep}

Pseudo-Anosov mapping classes have canonical representatives (up to topological conjugacy) which happen to be singular. The goal of this section is to describe how to symplectically smooth these singular representatives. We review the the structure of the singular representatives and then spend the rest of the section symplectically smoothing their singularities.

\begin{figure}
\includegraphics{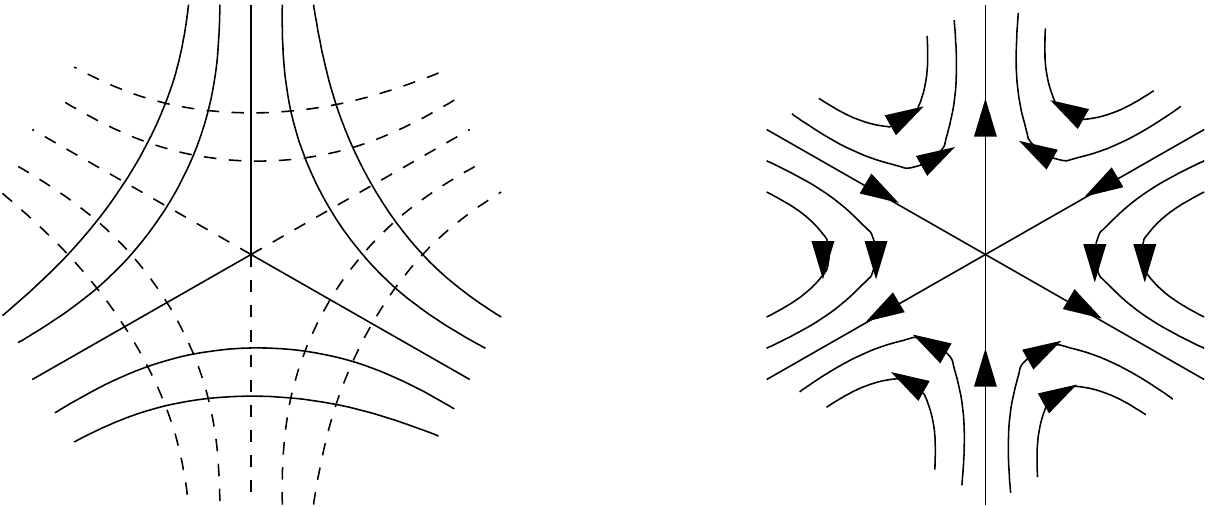}
\caption{\emph{Left}: The foliations $\mathcal{F}$ and $\mathcal{G}$ in the case $p = 3$. \emph{Right}: The singular Hamiltonian vector field whose time-$1$ flow is $\phi$ when $\rho = \id$ in the case $p = 3$.}
\label{flowanddotted}
\end{figure}

A singular measured foliation $\mathcal{F}$ on a surface $\Sigma$ is, away from a finite number of points, a smooth foliation with a transverse measure, that is, a measure of how far a curve is travelling in the direction transverse to the foliation. Away from the singular points the foliation has charts which take the leaves to horizontal lines in $\RR^2$ such that the transverse distance between leaves is simply the vertical distance in $\RR^2$. At a singular point, $\mathcal{F}$ has a number of prongs $p \geq 3$ as in the left of Figure \ref{flowanddotted} for the case $p = 3$ (more on this below).

A pseudo-Anosov map $\phi_{sing}$ (that is, a canonical singular representative) preserves two measured foliations $\mathcal{F}$ and $\mathcal{G}$ which are transverse everywhere except their singular points, which coincide and have the same number of prongs (as in the left of Figure \ref{flowanddotted}). The map $\phi_{sing}$ expands the transverse measure on $\mathcal{F}$ by a factor of $\lambda > 1$ and contracts the transverse measure on $\mathcal{G}$ by the same factor. See \cite{flp}, \cite{cb}, \cite{t}, \cite{m}, \cite{jg}, \cite{p} for more on pseudo-Anosov maps.

As in \cite[Expos\'e 9]{flp}, we can find a Markov partition for $\phi_{sing}$. That is, we can decompose our surface into a finite collection of homeomorphic images of rectangles, overlapping only on their boundaries, such that $\mathcal{F}$ and $\mathcal{G}$ pull back to the standard horizontal and vertical foliations. We equip these rectangles with the standard flat metric coming from the transverse measures. Upon gluing together the rectangles, we may obtain cone points of angle $p\pi$, which we call singular points with $p$ prongs.

We change coordinates by a homeomorphism isotopic to the identity such that the foliations are smooth away from the singular points. Furthermore we require that in a neighborhood of a singular point $x$ with $p$ prongs, we have a smooth chart such that one of our foliations $\mathcal{F}$ is given by the subbundle of the tangent space to $\Sigma$ on which the quadratic differential $z^{p-2} dz^2$ is positive real, and the other $\mathcal{G}$ by the subbundle on which it is negative real. See the left of Figure \ref{flowanddotted}.

We choose our symplectic form to be the area form associated with the flat-metric-with-cone-points described above. This is smooth even at the singularities. We must find a smooth symplectic perturbation of the map $\phi$ near the singularities. We describe the perturbation near a fixed singularity. Near a singularity that maps to a different singularity, one can use the work below under an identification of neighborhoods of the singularities.

\begin{figure}
\includegraphics{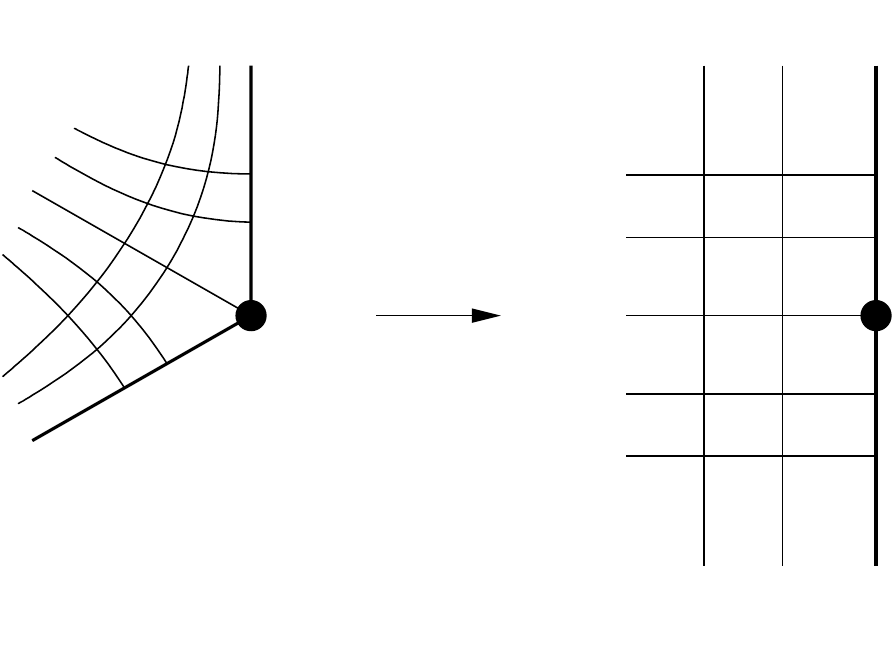}
\caption{Angular dilation of one of the regions transforms it to a standard rectangle with horizontal and vertical foliations.}
\label{dilation}
\end{figure}

If we divide up a neighborhood of this singular point by considering the components of the complement of the prongs of $\mathcal{F}$, each of these components is diffeomorphic to a rectangle on which the two foliations are the horizontal and vertical foliations. This diffeomorphism is given in radial coordinates by $(r,\theta)\mapsto (r,\frac{p}{2}\theta)$
 (see Figure \ref{dilation}). This map multiplies area by the constant multiple $\frac{p}{2}$ and thus a map is symplectic on one of the components if and only if the associated map is symplectic on the associated rectangle.

If $x$ is a fixed point, then $\phi$ cyclically permutes these rectangles. If we compose $\phi$ with an appropriate rotation $\rho$ (rotation by some multiple of $2\pi/p$) so that each rectangle is taken to itself, then the map $\rho\circ\phi = \phi\circ\rho$ on the closure of one of these rectangles is given in coordinates by the linear map $A_\lambda = \left(\begin{array} {cc}
\lambda & 0 \\
0 & \lambda^{-1} \end{array} \right)$, where the origin is the singularity. The map $A_\lambda$ is the time-$1$ flow of the Hamiltonian vector field associated to the Hamiltonian $H = \mu xy$, where $\mu = \ln \lambda$.

We can obtain a singular Hamiltonian $H_{sing}$ (singular at $x$ only) on an entire neighborhood of $x$ by piecing together the Hamiltonians on the various components. The Hamiltonians agree on the prongs, and the only issue is smoothness there. To see that we actually have a smooth Hamiltonian away from $x$, consider that if we use $\mathcal{G}$ instead of $\mathcal{F}$ and find Hamiltonians in the manner described above, they agree on the half-rectangles (given by the overlaps) with the Hamiltonians coming from $\mathcal{F}$ and they are smooth on the prongs associated to $\mathcal{F}$. Alternatively, we can simply note that $H_{sing} = \mu r^2\cos(p\theta)$ in polar coordinates, or $\mu \Re(z^p)/|z|^{p-2}$ in complex coordinates. See the right of Figure \ref{flowanddotted}.

Our strategy now is to modify $H_{sing}$ inside a small disk. We divide our task into two cases:
\bi
\ii $\rho \neq \id$ (the ``rotated'' case)
\ii $\rho = \id$ (the ``unrotated'' case)
\ei

\subsection*{Perturbing rotated singularities}
\label{rot}

Our map $\phi$ is given (in a neighborhood of $x$) by the time-$1$ flow of the singular Hamiltonian $H_{sing}$ composed with $\rho^{-1}$. We modify $\phi$ by forming $H_{sm} = f H_{sing}$, where $f$ is a smooth function whose zero set is a small ball and which is one outside a slightly larger one. Then we let $\phi_{sm}$ (in a neighborhood of $x$) be the time-$1$ flow of the Hamiltonian $H_{sm}$ composed with $\rho^{-1}$.

The union of all of the prongs is still preserved setwise by $\phi_{sm}$ and the rectangles continue to be permuted in the same manner as before. Thus the only possible fixed point is at $x$, where we have an elliptic fixed point, for $\phi_{sm}$ is given simply by $\rho^{-1}$ in a neighborhood of $x$.

\subsection*{Perturbing unrotated singularities}
\label{unrot}

Here $\phi$ is simply given (in a neighborhood of $x$) by the time-$1$ flow of $H_{sing}$. We must do something more elaborate than in the previous section in order to avoid degenerate fixed points.

We construct $H_{sm}$ in two steps. First we construct it on a small ball, where it has $p-1$ nondegenerate saddle points, and then we join it to $H_{sing}$ by an interpolation with no critical points.   See Figure \ref{smoothhsing} for the basic idea.

{\begin{figure}
\includegraphics{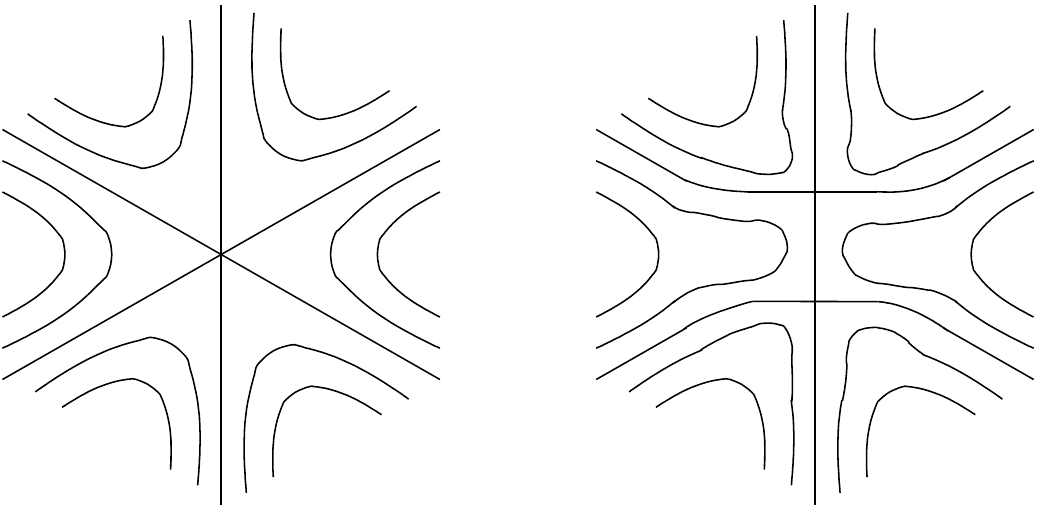}
\caption{Smoothing $H_{sing}$ in the case $p = 3$. On the left we have the level sets of $H_{sing}$ and on the right we have the level sets of the perturbed Hamiltonian $H_{sm}$ with $p - 1 = 2$ saddle points.}
\label{smoothhsing}
\end{figure}}

\begin{lemma}
\label{xsiny}
There is a smooth function $G$ on the closed disk $D = D(1)$ with $p-1$ nondegenerate saddle points and no other critical points such that on the boundary of the disk, $G|_{\del D}$ is Morse\footnote{We won't be doing Morse theory here; this is simply a convenient class of functions.} with $p$ minima and $p$ maxima. Additionally, consider the one-manifold $C_{max}$ consisting of the points where $G|_{\del D(r)}$ is maximum for $r$ near $1$. Likewise consider $C_{min}$. We require that the derivative of $G$ along $C_{max}$ (pointing outward) is positive, that the derivative along $C_{min}$ is negative, and that $C_{max}$ and $C_{min}$ are transverse to $\del D(r)$ where they are defined.
\end{lemma}

\proof
The function $x\sin(\pi y)$ has a saddle point at $(0,k)$ for $k \in \ZZ$. We consider a small rectangle $[-\epsilon,\epsilon]\times[0,p-2]$ with semicircular caps on each end. This is not smooth, but we may smooth it with a $C^1$-small perturbation near the points where it fails to be $C^2$. We claim the map on this domain $B$ has the requisite properties (we then choose an appropriate map from $D(1)$ to $B$).

On $B$, except near the ends, the maxima, minima, and zeroes near the boundary are (for an appropriate choice of the subdisks $D(r)$) horizontal lines and the map on the boundary has nondegenerate zeroes and extrema. 

The $C^1$-small perturbation from above can't produce extra critial points of $G|_{\del D(r)}$. Additionally, for $\epsilon$ small enough, $x\sin(\pi y)$ is a $C^2$-small perturbation of $\pi xy$, which has nondegenerate extrema and satisfies the derivative condition in the statement of the lemma, so $x\sin(\pi y)$ also has nondegenerate extrema. Finally, the map is seen to have the correct number of zeroes and extrema.
\qed

\begin{lemma}
\label{interp}
Let $G$ be as in Lemma \ref{xsiny}. Then there exist constants $0 < c_1 << 1$ and $c_2 >> 1$ such that if we define $H_{sm} = c_1 G(c_2 z)$ on a small disk of radius $1/c_2$, there exists an extension to a function $H_{sm}$ which equals $H_{sing}$ outside a somewhat larger disk, and which has no critical points in the intermediate region.
\end{lemma}

\proof
We proceed in three steps (we will rescale at the end):
\be

\ii Select $r_1$ slightly less than $1$ and consider the disk $D(r_1) \subset D(1)$. We modify $G$ so that the zeroes and extrema are in the correct places on $D(1) - D(r_1)$ (in the sense that their placement agrees with that of $H_{sing} = \mu r^2 \cos(p\theta)$).

Isotop the identity map on the disk (through maps which preserve $D(r)$ setwise) to one which takes the zeroes of $G|_{D(1)-D(r_1)}$ to the segments $\theta = \frac{k}{p} \pi/2$, the maxima of $G$ on the boundary of disks intermediate to $D_1$ and $D_2$ to the segments $\theta = \frac{2k}{p} \pi$, and the minima to the segments $\theta = \frac{2k+1}{p} \pi$. Note that these $\theta$ correspond with the placement of the zeroes, maxima, and minima of the map $H_{sing} = \mu r^2 \cos(p\theta)$.

Let $G_1$ be the composition of the inverse of the result of this isotopy with $G$. Note that $G_1$ also satisfies the conditions of Lemma \ref{xsiny}

\ii Now select $r_2$ such that $r_1 < r_2 < 1$ and consider the disk $D(r_2)$. We modify $G_1$ so that it is $f(r)\cos(p\theta)$ for some increasing function $f(r)$ on the region $D(1) - D(r_2)$.

To do this, select an increasing function $f(r)$ on $D(1)$ such that on $D(1) - D(r_1)$, $$\left|f(r)\cos(p\theta)\right| \geq K \left|G_1\right|$$ for some large constant $K$. 

Let $G_2$ be a modification of $G_1$ on $D(1) - D(r_1)$ which is increasing in absolute value (except where it is zero) on the rays $\theta = \textrm{const}$ as $r$ increases from $r_1$ to $1$ except possibly for a brief period initially (i.e. near $r_1$) and such that $G_2 = f(r)\cos(p\theta)$ on $D(1) - D(r_2)$. We claim that we do not create additional critical points. 
Because the maxima, minima, and zeroes have the desired behavior in Lemma \ref{xsiny}, there are only issues in between them. We can make 
the initial period where $G_2$ is not increasing short enough so as not to create new critical points of $G_2$ in this region.

\ii Extend $G_2$ to a function $G_3$ on $\RR^2$ by extending $f(r)$ to an increasing function which equals $C r^2$ for $r>>0$ and some large constant $C$. Finally, let $H_{sm}$ equal $c_1 G_3(c_2 z)$ for appropriate constants $c_1, c_2$ (defined on an appropriate region).
\ee

\qed

Thus our map $\phi_{sm}$ has $p-1$ positive hyperbolic fixed points, and no others, in the perturbed region. There are no fixed points coming from closed orbits because there are no components of level sets of $H_{sm}$ which are circles, for then there would necessarily be an extremum inside, and $H_{sm}$ has none.

\subsection{Floer homology of $\phi_{sm}$ for pseudo-Anosov maps}
\label{paproof}

We now consider the symplectic Floer homology of the smooth representative $\phi_{sm}$ we have just constructed. By Proposition \ref{notmono}, because the mapping class $h$ of $\phi_{sm}$ is pseudo-Anosov, $HF_*(\phi_{sm}) = HF_*(h)$, which in turn is equal to $HF_*(\phi)$ for any other $\phi$ in the same mapping class.

\begin{thm}
\label{nodiff}
In $CF_*(\phi_{sm})$, all differentials vanish. Thus $HF_*(h) = HF_*(\phi_{sm}) = CF_*(\phi_{sm})$.
\end{thm}

\proof
In \cite{bk} it is shown that no two fixed points (this includes both singular and nonsingular fixed points) of a singular, standard form pseudo-Anosov map $\phi_{sing}$ are Nielsen equivalent. If we can show that two fixed points of $\phi_{sm}$ are Nielsen equivalent if and only if they are both associated to the same singularity of $\phi_{sing}$, then we will be done, for a differential $u$ gives a Nielsen equivalence as in \S \ref{nclass}. Thus there are no differentials except between those associated to the same singularity. These, however, are all of the same index (mod 2)  by \S\ref{unrot}.

To see that it is indeed the case that two fixed points of $\phi_{sm}$ are Nielsen equivalent if and only if both are associated to the same singularity of $\phi_{sing}$, we argue as follows:

\be

\ii Whether two fixed points $x$ and $y$ of a map $\psi$ are Nielsen equivalent is unchanged by modifying a map inside a disk away from either fixed point (and thus any number of disks away from the fixed points), for any path passing through a disk is homotopic (rel boundary) to any other with the same endpoints.

\ii If two isolated fixed points $x$ and $y$ of a map $\psi$ are Nielsen inequivalent and we modify our map in a small disk near $x$ to a map $\psi'$, then any fixed point $x'$ (of the modified map) inside this disk is Nielsen inequivalent to $y$, for if we take a path $\gamma'$ exhibiting Nielsen equivalence of $x'$ and $y$ (under $\psi'$) and append a short path from $x$ to $x'$, then this new path $\gamma$ exhibits Nielsen equivalence of $x$ and $y$ (under our $\psi$).

To see this, note that the path following $\gamma$ from $x$ to $x'$ and then onward until it exits (for the last time) our small disk  at a point $x''$ (call this path $\alpha$) stays inside a slightly larger disk under $\psi$. Thus the homotopy from $x''$ to $\psi(x'') = \psi'(x'')$ given by the homotopy from $\gamma'$ to $\psi'(\gamma')$  extends to a homotopy rel $x$ from $\alpha$ to $\psi(\alpha)$. This implies that $\gamma$ is homotopic to $\psi(\gamma)$.

\ee

Repeated application of these two points gives our conclusion.

\qed

\begin{cor}
\label{pah1}
$H_1(\Sigma)$ acts as zero on $HF_*(h)$ for $h$ pseudo-Anosov.
\end{cor}

\proof
There are no flow-lines between generators of $CF_*(\phi_{sm})$, so the extrinsic or algebraic definition of the action of $H_1(\Sigma)$ as in Section \ref{product} is trivial. By Proposition \ref{productinvariant}, the same holds for $HF_*(h)$.
\qed

\section{Reducible maps}

\subsection{The space of $\eta$-weakly monotone maps}
\label{wmred}

In the reducible case, not every map is weakly monotone. We need to understand the structure of the space of $\eta$-weakly monotone maps, and in particular show that it is path connected, in order to prove invariance.\footnote{The results of this section go through for weakly monotone maps as well with similar arguments, which is sufficient for our needs, but involving Nielsen classes is not much more difficult.} We first recall the definition:

\begin{defn}
\label{etaweak}
The symplectomorphism $\phi: (\Sigma,\omega)\rightarrow(\Sigma,\omega)$ is \emph{weakly monotone} for a Nielsen class $\eta \in \pi_0(\Gamma(M_\phi))$  (or $\eta$-weakly monotone) if $\omega_\phi|_N = 0$, where $$N = \ker\left(c_\phi: H_1(\Gamma(M_\phi)_\eta) \rightarrow \RR\right),$$ where $\Gamma(M_\phi)_\eta$ denotes the $\eta$-component of $\Gamma(M_\phi)$.
\end{defn}

We would like to be able to define the space $\Symp_h^\eta(\Sigma,\omega)$ of maps $\phi$ in a mapping class $h \in \Gamma$ which are $\eta$-weakly monotone. To do this we must show that Nielsen classes are well-defined not only for a single map but also for an entire mapping class. Suppose we have an isotopy $\phi_t$ from $\phi = \phi_0$ to $\phi_1$. Thinking of a Nielsen class as a homotopy class of maps $f: [0,1] \rightarrow \Sigma$ with the restriction that $f(1) = \phi(f(0))$, we have a natural way to extend a given $f$ (up to reparametrization) to a representative of a Nielsen class for $\phi_t$ by extending the path $f$ from $f(0)$ to $\phi(f(0))$ to a path $f_t$ from $f(0)$ to $\phi_t(f(0))$ in the obvious manner. We must show that if $\phi_1 = \phi$ then the Nielsen class of $f_1$ is the same as the Nielsen class of $f$.

\begin{lemma}
\label{nielsenok}
Nielsen classes are well-defined on an entire mapping class if $\Sigma$ has negative Euler characteristic. That is, in the above setup, the Nielsen class of $f_1$ is the same as the Nielsen class of $f$.
\end{lemma}

\proof
We consider the universal cover $\HH$ of $\Sigma$. We have a map $\Phi : \Sigma\times S^1 \rightarrow \Sigma$ given by $(x,t)\mapsto \phi_t(x)$. This extends to a map $\widetilde\Phi: \HH\times [0,1] \subset \HH\times \RR \rightarrow \HH$ on the universal covers. If $\widetilde \Phi (\cdot, 0) = \widetilde \Phi (\cdot, 1)$, then the result follows. That is, we must show the map $\Phi \circ p : \HH\times S^1 \rightarrow \Sigma$ extends to a map to $\HH$ (where $p: \HH\times S^1 \rightarrow \Sigma\times S^1$ is the covering map).

Recall the homotopy lifting property, which states that this map extends if and only if $(\Phi\circ p)_* \pi_1(\HH\times S^1) = 0 \subset \pi_1(\Sigma)$. To see that this is the case, consider that $\pi_1(\Sigma\times S^1) = \pi_1(\Sigma) \times \ZZ$. The image that we're interested in is the image of the copy of $\ZZ$ under the map $\Phi_* : \pi_1(\Sigma) \times \ZZ \rightarrow \pi_1(\Sigma)$. Because $\phi$ is a diffeomorphism, $\Phi_* (\cdot, 0): \pi_1(\Sigma) \rightarrow \pi_1(\Sigma)$ is an isomorphism. Thus the image of the copy of $\ZZ$ is in the center of $\pi_1(\Sigma)$ because it commutes with every element of the image of $\Phi_* (\cdot, 0)$, which is all of $\pi_1(\Sigma)$.

Because $\chi(\Sigma) < 0$, the center of $\pi_1(\Sigma)$ is trivial,\footnote{To see this, note that if $\pi_1(\Sigma)$ has a center, then we get a copy of $\ZZ/k\times\ZZ/\ell$ acting freely and properly discontinuously by deck transformations (and, in particular, hyperbolic isometries) on $\HH$. The free hyperbolic isometries (i.e. the parabolic and hyperbolic ones) are infinite order, so $k = \ell = 0$. Thus the quotient is a $2$-manifold with fundamental group $\ZZ\times\ZZ$. The only one of these is the torus. This gives a hyperbolic structure on the torus, which violates Gauss-Bonnet.} and so the image of the copy of $\ZZ$ is indeed trivial. The result follows.
\qed

Thus we denote the space of maps in a mapping class $h \in \Gamma$ which are weakly monotone for $\eta$ by $\Symp_h^{\eta}(\Sigma,\omega)$. Similarly we use $\Symp_h^m(\Sigma,\omega)$ for the monotone symplectic maps in the mapping class $h$, and $\Symp_h(\Sigma,\omega)$ for all symplectic maps in the mapping class $h$. We note that if $\phi$ is monotone, then it is $\eta$-weakly monotone (for any $\eta$): the constant $k$ is simply the proportionality constant relating $\omega$ and $c_\phi = c_1(E)$.

\begin{prop}[{cf. \cite[Lemma 6 ff.]{s2}}]
\label{etahomotopic}
The inclusions $\Symp_h^m(\Sigma,\omega)\hookrightarrow\Symp_h^{\eta}(\Sigma,\omega)\hookrightarrow\Symp_h(\Sigma,\omega)$ are homotopy equivalences.
\end{prop}

\proof
For the duration of this proof, the $\Sigma$ and $\omega$ will be implied and we will omit them when they are not needed.

Fix $\phi \in \Symp_h^m$. We have the action $\Symp_0\times\Symp_h\rightarrow \Symp_h$. Using $\phi$, this gives a homeomorphism $\Symp_0\times\{\phi\}\rightarrow \Symp_h$. We ask which elements of $\Symp_0$ correspond under this homeomorphism to elements of $\Symp_h^\eta$.

Consider $\psi \in H_1(\Gamma(M_\phi)_\eta).$ We have a map $H_1(\Gamma(M_\phi)) \rightarrow H_2(M_\phi;\ZZ)$ given by taking the image in $H_2(M_\phi;\ZZ)$ of the fundamental class of the torus representing the element of $H_1(\Gamma(M_\phi))$. Denote the image of $\psi$ by $A_\psi \in H_2(M_\phi;\ZZ)$. Then if $\phi$ is $\eta$-weakly monotone, we have that $c_\phi(A_\psi) = 0$ implies $\omega_\phi(A_\psi) = 0$.

Suppose we perturb the map $\phi$ by an isotopy $\phi_t$ from $\phi = \phi_0$ to $\phi_1$. There is then a natural way to extend $\psi$ (up to reparametrization of $S^1$): $\psi$ is (as an element of $\Gamma(M_\phi)$) a path from $\gamma = \psi(0,\cdot)$ to $\phi^{-1}(\gamma)$. We let $\psi_t$ be the extension of this to a path from $\gamma$ to $\phi_t^{-1}(\gamma)$ in the obvious manner. We wish to understand how $\omega_{\phi_t}(A_{\psi_t})$ and $c_{\phi_t}(A_{\psi_t})$ vary with $t$.

The latter is constant: we are simply taking the Euler numbers of isomorphic bundles.

For the former, we briefly recall the Flux homomorphism (see \cite[\S 10.2]{ms} for details). This is a map $\textrm{Flux}: \widetilde{\Symp}_0(\Sigma,\omega)\rightarrow H^1(\Sigma;\RR)$ (where $\widetilde{\Symp}_0(\Sigma,\omega)$ is the universal cover of the identity component $\Symp_0(\Sigma,\omega)$, i.e. paths in $\Symp(\Sigma,\omega)$ starting at the identity up to homotopy). Its value, when paired with an element $\gamma\in H_1(\Sigma,\ZZ)$ represented by the smooth image of an $S^1$, is the area swept out by $\gamma$ under the path of symplectomorphisms. In our situation with $\chi(\Sigma) < 0$, we actually have $\textrm{Flux}: \Symp_0(\Sigma,\omega)\rightarrow H^1(\Sigma;\RR)$ because maps $S^1\times S^1 \rightarrow \Sigma$ are nulhomologous, and so the area swept out by these is zero.

Thus the difference between $\omega_{\phi_0}(A_{\psi_0})$ and $\omega_{\phi_1}(A_{\psi_1})$ is $$\textrm{Flux}(\{\phi_t^{-1}\circ\phi_0\}_{t\in[0,1]})[\phi_0^{-1}(\gamma)],$$ i.e. the area swept out (on $\Sigma$) by $\phi_t^{-1}(\gamma)$ as $t$ goes from $0$ to $1$. By the last line of the previous paragraph, we may simply write this as $\textrm{Flux}(\phi_1^{-1}\circ\phi_0)[\phi_0^{-1}(\gamma)]$.

Let $N$ be the subgroup of $H_1(\Sigma;\RR)$ generated by $[\phi^{-1}(\gamma)]$ for $\gamma$ the smooth image of an $S^1$ such that there exists a $\psi \in H_1(\Gamma(M_\phi)_\eta)$ with $\psi(0,\cdot) = \gamma$. Note that for each such $\gamma$, the class $A_\psi$ is well defined, for if we have two different $\psi$'s, putting them together we get a map $S^1\times S^1 \rightarrow \Sigma$, which then must be nulhomologous. Let $N_0 \subset N$ be generated by $[\phi^{-1}(\gamma)]$  for those $\gamma$ which additionally satisfy $c_\phi(A_\psi) = 0$.

We get that $\phi_1 \in \Symp_h^\eta$ if and only if $$\left.\textrm{Flux}(\phi_1^{-1}\circ\phi)\left(\phi^{-1}_*(\cdot)\right)\right|_{N_0}= 0.$$ 
That is, $\Symp_h^\eta$ is the kernel of the anti-homomorphism
$$\minCDarrowwidth20pt\begin{CD} \Symp_h @>\cdot^{-1} \circ \phi>> \Symp_0 @>\textrm{Flux}>> H^1(\Sigma;\RR) @>\cdot|_{N_0}>> \Hom(N_0,\RR)\end{CD}$$
This anti-homomorphism is surjective and continuous and thus all the fibers are homeomorphic (and homotopic in $\Symp_h$) and, seeing as the range is contractible, each fiber is also homotopic to all of $\Symp_h$. Thus $\Symp_h^\eta \simeq \Symp_h$ by the inclusion. As asserted in \cite[Lemma 6]{s2}, $\Symp_h^m \simeq \Symp_h$ by the inclusion (the proof is similar to what we've done here but is slightly simpler), and so $\Symp_h^m \simeq \Symp_h^\eta$ by the inclusion as well.
\qed

\begin{prop}
For any $\eta$-weakly-monotone symplectomorphism $\phi: (\Sigma,\omega)\rightarrow(\Sigma,\omega)$ in a mapping class $h$, $HF_*(\phi; \eta)$ is well-defined and $HF_*(\phi; \eta) = HF_*(h;\eta)$.
\end{prop}

\proof
We apply Proposition \ref{etahomotopic}: there is an isotopy $\phi_t$ with $\phi_0 = \phi$ and $\phi_1$ monotone such that $\phi_t$ is $\eta$-weakly monotone for all $t$. Now well-definedness follows by Condition \ref{ncond}. The final equality follows from Theorem \ref{invariance}.
\qed

\begin{cor}
\label{notmonored}
If $\phi: (\Sigma,\omega)\rightarrow (\Sigma,\omega)$ is a symplectomorphism in mapping class $h$ and is $\eta$-weakly monotone for all $\eta$, then $HF_*(\phi) = HF_*(h)$.\qed
\end{cor}

\subsection{Structure of reducible maps}
\label{red}

By Thurston's classification (see \cite{t} and \cite{flp}; also cf. \cite[Definition 8]{g}), in a reducible mapping class $h\in \Gamma$, there is a (not necessarily smooth)  map $\phi$ which satisfies the following:

\begin{defn}
\label{reduciblestandard}
A reducible map $\phi$ is in \emph{standard form} if there is a $\phi$-and-$\phi^{-1}$-invariant finite union of disjoint noncontractible (closed) annuli $N \subset \Sigma$ such that:
\be
\ii \label{twist} For $A$ a component of $N$and $\ell$ the smallest positive integer such that $\phi^\ell$ maps $A$ to itself, the map $\phi^\ell|_A$ is either a \emph{twist map} or a \emph{flip-twist map}. That is, with respect to coordinates $(q,p)\in[0,1]\times S^1$, we have one of \beqa(q,p)&\mapsto& (q,p-f(q)) \qquad \textrm{(twist map)} \\ (q,p)&\mapsto& (1-q,-p+f(q)) \qquad \textrm{(flip-twist map)},\eeqa
where $f : [0,1]\rightarrow \RR$ is a strictly monotonic smooth map. We call the (flip-)twist map \emph{positive} or \emph{negative} if $f$ is increasing or decreasing, respectively. Note that these maps are area-preserving.
\ii Let $A$ and $\ell$ be as in (\ref{twist}). If $\ell = 1$ and $\phi|_N$ is a twist map, then $\Im(f)\subset[0,1]$. That is, $\phi|_{\textrm{int}(A)}$ has no fixed points. (If we want to twist multiple times, we separate the twisting region into parallel annuli separated by regions on which the map is the identity.) We further require that parallel twisting regions twist in the same direction.
\ii For $S$ a component of $\Sigma \backslash N$ and $\ell$ the smallest integer such that $\phi^\ell$ maps $A$ to itself, the map $\phi^\ell|_S$ is area-preserving and is either \emph{periodic} (i.e. some power is the identity) or \emph{pseudo-Anosov}. If periodic, we require it to have the form of an isometry of a hyperbolic surface with geodesic boundary. We define pseudo-Anosov on a surface with boundary and standard form for such below. If a periodic component $S$ has $\ell = 1$ and $\phi|_S = id_S$ we will furthermore call it \emph{fixed}.
\ee
\end{defn}

\begin{remark}
\label{puncspecial}
When a pseudo-Anosov component abuts a fixed component with no twist region in between, we'll need to deal with this case differently throughout. This is the major interaction between components. (There is also the influence of the twisting direction in neighboring twist regions on fixed components.) When two pseudo-Anosov components meet, we perturb both as below and there is no interaction.
\end{remark}

A model for a pseudo-Anosov map on a surface with boundary is given in \cite[\S2.1]{jg}. The idea is to define a pseudo-Anosov map on a surface with punctures and then ``blow up'' the punctures to recover boundary components. To define a pseudo-Anosov map on a surface with punctures, we simply consider the punctures as (yet more) singularities and let our foliations $\mathcal{F}$ and $\mathcal{G}$ be given by the subbundles of the tangent space of $\Sigma$ on which the quadratic differential $z^{p-2}dz^2$ is positive (resp. negative) real, but instead of requiring $p \geq 3$, we only require $p \geq 1$. The case $p = 2$ corresponds to a smooth fixed point (which we have punctured at the fixed point). As in Section \ref{symprep}, the map $\phi$ is given in polar coordinates in a neighborhood of the puncture by the time-$1$ flow of the (singular at the puncture) Hamiltonian $H_{punc} = \mu r^2 \cos(p\theta)$. 

We must now blow up each puncture to a boundary circle such that the map on the boundary $S^1$ is a rotation by some angle $\alpha$ (to match up with the maps in item (\ref{twist}) in Definition \ref{reduciblestandard}). This can be done in various ways (see \cite[\S2.1]{jg} for example), but we must do it in such a way that the resulting map is area-preserving and all fixed points are non-degenerate. Note that we are free to choose $\alpha$ as the difference can be made up by twisting (or by undoing twisting).

As in Section \ref{symprep}, we again have our rotation $\rho$ (which is rotation by some multiple of $2\pi/p$) and this again gives us two cases:
\bi
\ii $\rho \neq \id$ (the ``rotated'' case)
\ii $\rho = \id$ (the ``unrotated'' case)
\ei

\subsection*{Perturbing rotated punctures}
\label{rotpunc}

In the rotated case, we choose $\alpha$ to be the same as the rotation angle for $\rho^{-1}$. The same argument as in Section \ref{rot} allows us to modify the map so that it has no fixed points except inside a disk which is rotated with angle $\alpha$. Excising a sub-disk, we have our boundary component with rotation angle $\alpha$.

\subsection*{Perturbing unrotated punctures}
\label{unrotpunc}

The unrotated case is again more complicated. In order for the fixed points to be non-degenerate, we'll need the rotation angle $\alpha$ to be a non-trivial one. The desired result will be a Hamiltonian whose flow is a rotation inside a small ball. We'll find one which has $p$ hyperbolic fixed points around the rotating region, and agrees with the flow of $H_{punc}$ outside a slightly larger ball. See Figure \ref{smoothpunc} for the basic idea.

{\begin{figure}
\includegraphics{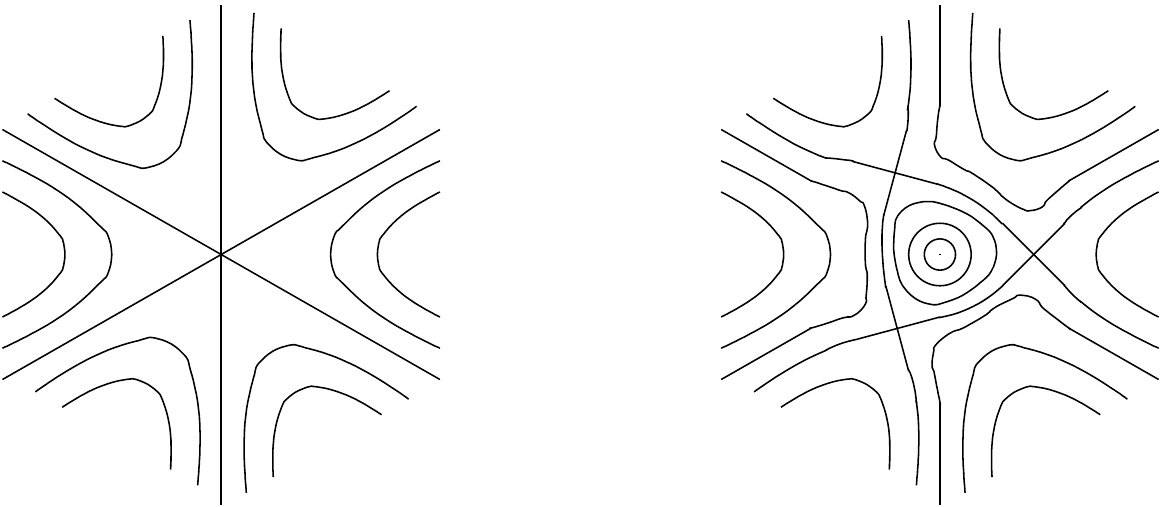}
\caption{Modifying $H_{punc}$ in the case $p = 3$. On the left we have the level sets of $H_{punc}$ and on the right we have the level sets of the perturbed Hamiltonian $H_{sm}$ whose flow is rotation inside a small ball (with an elliptic fixed point at the center, but we're excising the ball so it won't count) and has $p = 3$ hyperbolic fixed points surrounding the rotating region.}
\label{smoothpunc}
\end{figure}}

We proceed in four steps. First we add a positive bump of the form $f(r,\theta) = Ae^{-kr^2}$ to $H_{punc}$ to produce the $p$ equally spaced saddle points. Then we smoothly cut off the tail of $f$ without modifying the critical points of $H_{punc} + f$. Next we modify $H_{punc} + f$ near the origin so that its flow rotates at constant angular velocity. Finally we ensure there are no fixed points of the time-$1$ flow of the Hamiltonian vector field except for those corresponding to our $p$ critical points.

\emph{Step 1}: We consider $H_1 = H_{punc} + f(r,\theta) = \mu r^2 \cos(p\theta) + Ae^{-kr^2}$ with $A$ and $k$ positive constants to be chosen later. We're interested in the critical points of this function (away from the origin, where it's singular), so we calculate its partial derivatives:

\beqa
\frac{\del H_1}{\del r} &=& 2\mu r \cos(p\theta) - 2krAe^{-kr^2} \\
\frac{\del H_1}{\del \theta} &=& -p\mu r^2\sin(p\theta)
\eeqa

Thus at the critical points of $H_1$ we must have $\sin(p\theta) = 0$. In this case, $\cos(p\theta) = \pm 1$. If it's $-1$, $\frac{\del H_1}{\del r}$ cannot be zero (because both terms are negative). When it's $1$, $\frac{\del H_1}{\del r} = 2\mu r - 2krA e^{-kr^2}$. For $r > 0$, this is zero when $e^{kr^2} = \frac{kA}{\mu}$, i.e. when $r = r_c := \sqrt{\frac{1}{k}\ln(\frac{kA}{\mu})}$. We impose the restriction that $kA > \mu$. Note that by making $k$ large, we can make $r_c$ arbitrarily small.

In summary, we get critical points at one value of $r$ at the values of $\theta$ when $\cos(p\theta) = 1$, that is, for $p$ values of $\theta$. To check that these are all saddle points, we compute the Hessian at these points:
\beqa
\frac{\del^2 H_1}{\del r^2} &=& 2\mu + 2\mu k r^3 > 0\\
\frac{\del^2 H_1}{\del \theta^2} &=& -p^2\mu r^2 < 0\\
\frac{\del^2 H_1}{\del r \del \theta} &=& 0.
\eeqa
Thus the critical points are all nondegenerate and of index one, so we have $p$ saddle points.


\emph{Step 2}: Keeping $f$ solely a function of $r$, we cut it off smoothly starting at some point past $r_c$ to give a Hamiltonian $H_2$ which agrees with $H_{punc}$ outside a ball. As long as we keep $\frac{\del f}{\del r} > -2\mu r$, we create no new critical points.

Note that $f(r_c) = \frac{\mu}{k}$. Keeping $\frac{\del f}{\del r}$ near $-\mu r_c$ (which, using e.g. $A = \frac{e\mu}{k}$, is $\frac{1}{\sqrt{k}}$), we can bring $f$ to zero in a radial distance of a constant times $\frac{1}{\sqrt{k}}$; i.e. for $k$ large we can make $H_2$ agree with $H_{punc}$ outside an arbitrarily small ball.

\emph{Step 3}: Now we modify $H_2$ near the origin to give us $H_3$ which is $-Br^2 + C$ near the origin (for $B$ positive), which corresponds to the Hamiltonian flow rotating at a constant angular rate. Since $\frac{\del H_2}{\del r}$ is negative for $r < r_c$, we can patch together $-Br^2 + C$ near the origin with $H_2$ outside a small ball (of radius less than $r_c$) in a radially symmetric manner to get $H_3$ such that $\frac{\del H_3}{\del r}$ is negative for $r<r_c$ (we do this by choosing $C$ sufficiently large).

\emph{Step 4}: Finally, to ensure no fixed points of the time-$1$ flow of $H_{sm}$, we let $H_{sm}$ be $H_3$ multiplied by a radially symmetric function which is $\epsilon$ for $r < R$ (for $\epsilon$ sufficiently small that the only fixed points of the time-$1$ flow inside radius $R$ are the critical points and for $R$ large enough that $H_3$ agrees with $H_{punc}$ for $r > R$) and $1$ for $r > 2R$ (we may have to make $r_c$ smaller so we can fit a ball of radius $2R$, but this is no problem). This creates no new fixed points in the region $R < r < 2R$ because $H_3$ and $\frac{\del H_3}{\del r}$ have the same sign there. Now there are no fixed points of $\phi_{sm}$, the time-$1$ flow of the Hamiltonian vector field of $H_{sm}$, except for the $p$ critical points of $H_{sm}$ because outside radius $R$ there are no compact flow lines.

Thus our map $\phi_{sm}$ has $p$ positive hyperbolic fixed points, and no others, in the perturbed region. Note that everything in this section works in exactly the same manner if we add on a small negative bump instead of a small positive bump, which we may do in order to twist in the opposite direction.

\subsection{Nielsen classes for reducible maps}
\label{nielsenred}

We now study the Nielsen classes of fixed points of standard form reducible maps. 
Let us briefly describe the fixed points of our standard form reducible maps. We have:

\bi
\ii (Type Ia) The entire component of fixed components $S$ of $\Sigma\backslash N$ with $\chi(S) < 0$.
\ii (Type Ib) The entire component of fixed components $S$ of $\Sigma\backslash N$ with $\chi(S) = 0$. These are annuli and only occur when we have multiple parallel Dehn twists.
\ii (Type IIa) Fixed points $x$ of periodic components $S$ of $\Sigma\backslash N$ with $\chi(S) < 0$ which are setwise fixed by $\phi$. (These can be understood by considering the map on $S$ as a hyperbolic isometry, from which we see that $x$ must be an elliptic fixed point).
\ii (Type IIb) Fixed points $x$ of flip-twist regions. (These are elliptic. Note that each flip-twist region has two fixed points.)
\ii (Type III) Fixed points $x$ of pseudo-Anosov components $S$ of $\Sigma\backslash N$ which are setwise fixed by $\phi$. These come in 4 types (note that there are no fixed points associated to a rotated puncture):
\bi
\ii (Type IIIa) Fixed points which are not associated with any singularity or puncture (i.e. boundary component) of the pre-smoothed map. These may be positive or negative hyperbolic.
\ii (Type IIIb-$p$) Fixed points which come from an unrotated singular point with $p$ prongs. There are $p-1$ of these for each such, all positive hyperbolic.
\ii (Type IIIc) Fixed points which come from a rotated singular point. There is one for each such and it is elliptic.
\ii (Type IIId-$p$) Fixed points which come from an unrotated puncture (boundary component) with $p$ prongs. There are $p$ for each such, all positive hyperbolic.
\ei
\ei

In the smooth case (as opposed to the area-preserving case, which we're dealing with), Type Ib fixed points can be isotoped away. To see why, we introduce the concept of a \emph{multiple twist region} which is a (maximal) annulus which is the union of a collection of twist annuli and fixed annuli (i.e. those of type Ib); note that every fixed annulus is between two twist annuli. In other words, this is a region in which we're performing multiple parallel Dehn twists. If we're allowed non-area-preserving maps, we can get rid of all fixed points in such a region by performing an isotopy which moves every interior point closer to one boundary component.

With this setup, the Nielsen classes of fixed points of maps which are in standard form except that:
\bi
\ii All Type Ib fixed points have been eliminated (in the above manner)
\ii Fixed points of Type IIIb associated to the same singular point are left as one fixed point
\ii Fixed points of Type IIId associated to the same puncture are left as a full circle (i.e. the blown-up puncture)
\ei
have been studied by Jiang and Guo \cite{jg}. They show that in this case, Nielsen classes of fixed points are all connected --- that is, two fixed points are in the same Nielsen class if and only if there is a path between them through fixed points. That is, we have a separate Nielsen class for every component of Type Ia, for every single fixed point of Type IIa, IIb, IIIa, or IIIc, and for every unrotated singular point of the pre-smoothed map for Type IIIb (i.e. the collection of fixed points associated to a single unrotated singular point are all in the same Nielsen class).

Type IIId is special, as mentioned in Remark \ref{puncspecial}. If a pseudo-Anosov component meets a fixed component, then notice that the full circle of fixed points of Type IIId (in the above non-area-preserving model) coincides with fixed points on the boundary of the Type Ia fixed points (note that a pseudo-Anosov component never meets a Type Ib component, as these only occur between twist regions). In this case there will be some interaction, and, while we will perturb the pseudo-Anosov side as in Section \ref{unrotpunc} and the fixed side with a corresponding small rotation near the boundary, the resulting $p$ fixed points will be in the same Nielsen class as the Type Ia fixed points in the fixed component. If two pseudo-Anosov components meet, the two full circles coincide and, after perturbing as in Section \ref{unrotpunc}, we will have $p + q$ positive hyperbolic fixed points all in the same Nielsen class, where $p$ is the number of prongs on one side and $q$ is the number of prongs on the other side. Note that in this case we will have to twist in a coherent manner, but this is not an issue. In all other cases, we have a separate Nielsen class associated to the boundary component of the pseudo-Anosov  component containing $p$ positive hyperbolic fixed points.

Jiang and Guo's argument actually extends to imply that, if we don't eliminate Type Ib fixed points (which we won't be able to do in an area-preserving manner), there is a separate Nielsen class (separate from all the others just mentioned) for each component of Type Ib. We briefly explain this.\footnote{Gautschi explains how to do this as well \cite[Proposition 20]{g}, but his argument is specific to the case in which there are no pseudo-Anosov components.} We note that the argument will crucially use the fact that in multiple twist region, all twisting happens in the same direction (i.e. all parallel Dehn twists have the same sign).

Recall that two fixed points $x$ and $y$ are Nielsen equivalent if there exists a path $\gamma(t)$ with $\gamma(0) = x$, $\gamma(1) = y$, and $\phi(\gamma(t))$ is homotopic rel boundary to $\gamma(t)$. This notion extends to $\phi$-invariant sets (i.e. setwise fixed). We say that two $\phi$-invariant sets $A$ and $B$ (either of which may be a single fixed point) are $\phi$-related if there exists a path $\gamma : ([0,1],0,1) \rightarrow (\Sigma,A,B)$ such that $\phi(\gamma) \simeq \gamma$ through maps $([0,1],0,1) \rightarrow (\Sigma,A,B)$. We say that such a $\gamma$ is a $\phi$-relation between $A$ and $B$.

For a map $\phi$ in standard form, we consider a collection of $\phi$-invariant reducing curves $\{C_i\}$ given by $\del N$ for $N$ our collection of annuli, excepting that when we have parallel annuli, we only use the two outermost $\phi$-invariant curves; that is, we consider a multiple twist region as a single annulus. The following is (part of what is) proved in \cite[\S3.2-3]{jg}:

\begin{lemma}[{\cite[\S3.2-3]{jg}}]
\label{jg}
Suppose $x$ and $y$ are two fixed points of a map $\phi$ which is in standard form which are Nielsen equivalent. Let $\tilde\gamma$ be a $\phi$-relation between $x$ and $y$. Homotop $\tilde\gamma$ to a $\phi$-relation $\gamma$ with minimal (geometric) intersection number with the $C_i$ in its homotopy class. Break this path into segments at its crossings with the $C_i$. Let $\gamma_j$ be the $j^{\textrm{th}}$ segment and let $S_j$ be the (closure of the) component of $\Sigma \backslash \{C_i\}$ on which $\gamma_j$ lives ($S_j$ may be an annulus). If $\gamma_j$ is an interior segment, let $A_j$ and $B_j$ be the two reducing curves (i.e. elements of $\{C_i\}$) the segment $\gamma_j$ intersects. If $\gamma_j$ is an initial segment, instead let $A_j = x$ and if $\gamma_j$ is a terminal segment, let $B_j = y$. Then $\gamma_j$ is a $\phi$-relation between $A_j$ and $B_j$ on the subsurface $S_j$. Furthermore, $A_j$ and $B_j$ must be pointwise fixed. Furthermore, we may homotop $\gamma_j$ inside the fixed point locus unless we're dealing with a fixed point (i.e. $A_j = x$ or $B_j = y$) and the fixed point is in a multiple twist region. (The same is true here, see the Corollary below.)
\qed
\end{lemma}

\begin{cor}
\label{multitwist}
For a map $\phi$ in standard form, each component of Type Ib fixed points is in its own Nielsen class, distinct from the Nielsen class of any other fixed point of $\phi$. Furthermore, if $\gamma$ is a $\phi$-relation between two fixed points lying inside a multiple twist region, we may homotop it to lie inside the fixed point locus.
\end{cor}

\proof
By Lemma \ref{jg}, we need only show that on the closed annulus on which we've performed multiple parallel Dehn twists all of the same sign, Nielsen classes of fixed points are connected (i.e. they are simply the fixed sub-annuli and, if fixed, boundary circles). Note that we allow the boundary to twist or remain fixed.

Consider a path $\gamma$ between two fixed points $x$ and $y$ in this situation.  Consider $1$-chains with fixed endpoints $x$ and $y$ up to homology. This is an affine space over the first homology of the annulus, which is $\ZZ$. The difference between $\phi(\gamma)$ and $\gamma$ in terms of this relative homology is the element of the homology of the annulus corresponding to the number of twisting regions between $x$ and $y$, with sign. Because all twisting regions have the same sign, this number is not zero unless $x$ and $y$ are in the same component of the fixed point locus of $\phi$, and so they cannot be homotopic unless $x$ and $y$ are in the same component.

If $x$ and $y$ are in the same component of the fixed point locus, then we can easily homotop our $\phi$-relation $\gamma$ to lie inside the component (because the annulus deformation retracts onto the circle), proving the second statement as well.
\qed

\subsection{Weak monotonicity of standard form reducible maps}
\label{redwkmono}

We now show that $\phi = \phi_{sm}$ in standard form is weakly monotone.
\begin{lemma}
\label{noperiods}
Let $\phi: \Sigma \rightarrow \Sigma$ be in standard form. Then $\omega_{\phi}: H_1(\Gamma(M_\phi)) \rightarrow \RR$ is the zero map.
\end{lemma}

\proof
Consider a section $\psi$ of $$\begin{CD}
\Sigma @>>> S^1\times M_\phi \\
@. @VVV \\
@. S^1\times S^1 \\
\end{CD}$$ such that $\psi(0,\cdot)$ is in Nielsen class $\eta$ and $[\psi] \in H_1(\Gamma(M_\phi)_\eta)$ is the class we wish to test $\omega_{\phi}$ on.

As in Lemma \ref{pi1is0}, let $\alpha_t(s) = \psi(s,t)$. Each $\alpha_t$ should be thought of as a closed curve on $\Sigma$ by identifying all fibers of $M_\phi = (\RR\times\Sigma)/\hspace{-5pt}\sim\ \longrightarrow S^1$ with the fiber over $t = 0$ by projecting $\RR\times\Sigma$ to $\Sigma$ (with $t = 0$ and $t = 1$ thought of as separate times). Then $\phi(\alpha_0) = \alpha_1$ is homotopic to $\alpha_0$ through $\alpha_t$. Let $A_\psi$ denote the $2$-chain with boundary given by this homotopy. Note that this is equivalent to giving the homology class $[\psi(S^1\times S^1)]$ (the two are related by appending the tube $S^1\times \alpha_0$ to the $2$-chain $A_\psi$, thought of as living in the fiber $\Sigma$ over zero). We have $\omega_\phi(A_\psi) = \int_{S^1\times S^1} \psi^* \omega_\phi$.

We desire to show that
$\omega_\phi(A_\psi) = 0$.
Note that continuously varying the homotopy $\alpha_t$ has no effect on $\omega_\phi(A_\psi)$.
Nor does continuously varying $\alpha_0$ (which requires that we vary $\alpha_1 = \phi(\alpha_0)$ appropriately) as this does not change the cohomology class $[\psi(S^1\times S^1)]$.

\emph{Claim}: There is a choice of $\psi$ such that $A_\psi$ either lies entirely inside a single pseudo-Anosov component or avoids all pseudo-Anosov components.

Given this, we have two cases.

\emph{Case 1}: $A_\psi$ lies entirely in $S$, a pseudo-Anosov component. We show that $\omega_\phi(A_\psi) = 0$ (i.e. the constant $k$ is zero). In this case, by Lemma \ref{nohom}, we must have $\alpha_0$ either nulhomologous or boundary parallel.\footnote{Technically, that result was for pseudo-Anosov maps on closed surfaces. Jiang and Guo give an argument \cite[Proof of Lemma 2.2]{jg} which includes pseudo-Anosov maps with boundary.} If it's nulhomologous, we may argue as in Lemma \ref{pi1is0} and conclude that $\omega_\phi(A_\psi) = 0$.

If $\alpha_0$ is boundary parallel, we homotop it near the boundary component so that the map there is given by the local model of Section \ref{unrotpunc} (see also Figure \ref{smoothpunc}). This map is then Hamiltonian, and so, choosing our homotopy to be one also lying in this region (which, as noted above, does not change the class of $A_\psi$) we see that $\omega_\phi(A_\psi) = 0$.

\emph{Case 2}: $A_\psi$ avoids all pseudo-Anosov components. In this case, we may as well be considering a map in standard form which has no pseudo-Anosov components (for example, by replacing pseudo-Anosov components with periodic components with the same boundary behavior). This is the situation Gautschi considers, and he shows in \cite[Proposition 13]{g} that $\omega_\phi(A_\psi) = 0$ for all $\psi$ in this case.

We now prove the above claim. We consider two cases.

\textbf{Case A: $\phi$ has a fixed point in class $\eta$}

We show that we may assume $A_\psi$ lies entirely inside the component $S$ of $\Sigma\backslash N$ in which $x$ lies. This requires showing both that we can restrict to $\alpha_0 \subset S$ (which implies $\alpha_1 = \phi(\alpha_0) \subset S$) and that we can bring the homotopy $\alpha_t$ inside $S$.

Notice that $\alpha_0$ is a $\phi$-relation between $x$ and itself. As in the statement of Lemma \ref{jg}, 
we homotop $\alpha_0$ so that it has minimal geometric intersection with a reducing set of curves. Then by Lemma \ref{jg} and the fact that the collection of fixed points in a single Nielsen class is a subset of a single component of $\Sigma\backslash N$ (see Section \ref{nielsenred}), $\alpha_0$ must lie in $S$.

To see that we can bring the homotopy $\alpha_t$ inside $S$, notice that the inclusion $S\hookrightarrow \Sigma$ is $\pi_1$-injective. Thus there exists a possibly different homotopy $\alpha'_t$ between $\alpha_0$ and $\alpha_1$ which lies in $S$. Piecing these together we get a map $S^1\times S^1 \rightarrow \Sigma$ which is necessarily nulhomologous because $\chi(\Sigma) < 0$. Thus the class $A_\psi$ is the same for our new homotopy. 

\textbf{Case B: $\phi$ has no fixed point in class $\eta$}

We argue in 3 steps that $A_\psi$ may be homotoped either entirely inside a pseudo-Anosov component or to avoid all pseudo-Anosov components:

\emph{Step 1}: Let $\{C_i\}$ be the subset of the reducing curves given by the boundaries of the pseudo-Anosov components (including those abutting Type Ia components). Inspired by Jiang and Guo \cite[\S3.2-3]{jg} (cf. Lemma \ref{jg}), we homotop $\alpha_0$ so it has minimal geometric intersection with the collection of the $C_i$ and work in the universal cover $\HH$ of $\Sigma$. Each $C_i$ is noncontractible, so each lift in $\HH$ is a copy of $\RR$, and each component $S_k$ of $\Sigma \backslash \coprod_i C_i$ is $\pi_1$-injective, so each lift of $S_k$ is also a universal cover of $S_k$. Furthermore, each lift of a $C_i$ separates $\HH$ into two components. Thus $\widetilde\alpha_0$ (a chosen lift of $\alpha_0$) intersects each lift of $C_i$ at most once (else we could perform a finger move to reduce intersections).

We lift $\phi$ to $\widetilde\phi: \HH \rightarrow \HH$ and note that $\widetilde\phi$ maps lifts of $C_i$ to lifts of some $C_j$ (for each $i$ and some $j$). Thus $\widetilde\phi(\widetilde\alpha_0)$ also has minimal intersection with the lifts of the $C_i$.

\emph{Step 2}: Choose some basepoint $p$ on $\alpha_0$ and let $\widetilde p$ be a lift compatible with our chosen lift of $\widetilde \alpha_0$. Consider $\beta(t) = \psi(0,t)$, a curve from $p$ to $\phi(p)$. Choose a lift $\widetilde\beta$ of $\beta$ through $\widetilde p$. Let $\widetilde\alpha$ be a segment inside $\widetilde\alpha_0$ from $\widetilde p$ to $\widetilde p'$ (another lift of $p$) which is injective on its interior (i.e. its projection wraps around $\alpha$ once).

Thus we have a rectangle in the universal cover with corners $\widetilde p$, $\widetilde p'$, $\widetilde \phi(\widetilde p)$, and $\widetilde\phi(\widetilde p')$ and edges $\widetilde \alpha$, $\widetilde\phi(\widetilde\alpha)$, $\widetilde \beta$, and $\widetilde \beta'$ (the lift of $\widetilde \beta$ through $\widetilde p'$). We call its interior $\widetilde \psi$.

We will be homotoping these objects into easier to work with positions but will not change their names as we do so: homotop $\beta$ (and thus also $\widetilde\beta$ and $\widetilde\beta'$) and thus also $\alpha_0$, so it has minimal intersection with the $C_i$, so the lifts each intersect each lift of a $C_i$ at most once. We can do this without disturbing this property for $\alpha_0$.

\emph{Step 3}: Homotop $\widetilde\psi$ rel boundary so that it is transverse to each lift of a $C_i$. Let $\widetilde C$ be the intersection of $\widetilde \psi$ with the lifts of the $C_i$. This is a compact one-manifold with boundary in the rectangle $\widetilde\psi$.

Note that $\widetilde C$ has no circle components because all of the lifts of the $C_i$ are copies of $\RR$. Also note that no component of $\widetilde C$ has both boundary components on the same edge of the rectangle $\widetilde \psi$ because each of the edges intersects each lift of a $C_i$ at most once.

\emph{Step 4}: Each of the edges of $\psi$ has a certain number of boundary points of $\widetilde C$ on it. Note that opposite edges have the same number of such (which are in corresponding positions) because one pair is a lift of the same path $\beta$ on $\Sigma$, and the other pair are lifts of curves $\alpha_0$ and $\phi(\alpha_0)$ and $\phi$ preserves the collection $\{C_i\}$.

Using this and Step 3, we have three cases:

\underline{Case 1}: There are no boundary points of $\widetilde C$ on $\widetilde \alpha$ (and thus also on $\widetilde\phi(\widetilde \alpha)$).

Thus $\widetilde C$ consists of some number of (non-intersecting) paths from $\widetilde\beta$ to $\widetilde\beta'$ which each go from a point on $\widetilde\beta$ to the corresponding point on $\widetilde\beta'$. If we project to $\Sigma$, what we see is the annulus given by $A_\psi$ with boundary $\alpha_0$ and $\phi(\alpha_0) = \alpha_1$ with some number of intermediate $C_i$'s. If there is at least one intermediate $C_i$, then $\alpha_0$ is homotopic to a $C_i$. If we homotop it so that its image is $C_i$, we find that $\omega_\phi(A_\psi) = 0$, as desired.

If there are no intermediate $C_i$'s, then $A_\psi$ is contained in one component of $\Sigma\backslash \coprod_i C_i$, as desired.

\underline{Case 2}: There are no boundary points of $\widetilde C$ on $\widetilde \beta$ (and thus also on $\widetilde \beta'$).

In this case, $\widetilde C$ consists of some number of (non-intersecting) paths from $\widetilde \alpha$ to $\widetilde\phi(\widetilde\alpha)$ which each go from a point on $\widetilde\alpha$ to a point on $\widetilde\phi(\widetilde\alpha)$. If we project to $\Sigma$, what we see is the annulus given by $A_\psi$ with boundary $\alpha_0$ and $\phi(\alpha_0) = \alpha_1$ with some number of $C_i$'s (which are non-intersecting) running from a point on $\alpha_0$ to the corresponding point on $\alpha_1$. If there are no such $C_i$'s, then $A_\psi$ is contained in one component of $\Sigma\backslash \coprod_i C_i$, as desired. Similarly if there is only one $C_i$.

Suppose there are at least two such $C_i$. We desire to show that $A_\psi$ avoids pseudo-Anosov components. Suppose that it does not. Then there is some rectangle in $A_\psi$ with boundary edges given by portions of some $C_i$ and some $C_j$ as well as portions of $\alpha_0$ and $\alpha_1$ which is contained in a pseudo-Anosov component. Note that $\phi$ applied to this portion of $\alpha_0$ is the portion of $\alpha_1$.

Work only in the pseudo-Anosov component as a surface with boundary. If $C_i$ (resp. $C_j$) corresponds to a boundary component with $p > 1$ prongs, consider a surface in which we've collapsed $C_i$ (resp. $C_j$) and have a pseudo-Anosov map on this surface with fewer boundary components. Near any boundary component with one prong, we have a (hyperbolic) fixed point nearby. After collapsing any with $p > 1$, we have a fixed point of some kind nearby (we may perturb to our standard map if desired; there will still be a fixed point of some kind nearby). By the argument in the proof of Theorem \ref{nodiff}, we may modify $\alpha_0$ slightly to exhibit a Nielsen equivalence between a fixed point near $C_i$ and a fixed point near $C_j$ (recall that $\alpha_1 = \phi(\alpha_0)$). Thus by the results of Section \ref{nielsenred} (or by \cite{bk}, as in the proof of Theorem \ref{nodiff}, if our surface is now closed), $C_i = C_j$. Now if $C_i$ has $p > 1$ prongs so that we've collapsed it, then we have a curve $\gamma$ on the surface such that $\phi(\gamma) = \gamma$ (using the induced map $\phi$). Thus $\gamma$ is nulhomotopic by \cite[Lemma 2.2]{jg}. Similarly, if $C_i$ is a $p = 1$ boundary component and again by \cite[Lemma 2.2]{jg} we find that it is contractible. But this is impossible ---  we can reduce intersections of $\alpha_0$ with $C_i = C_j$ by a finger move. Thus $A_\psi$ avoids pseudo-Anosov components, as desired.

\underline{Case 3}: There are boundary points of $\widetilde C$ on both $\widetilde \alpha$ and $\widetilde \beta$.

In this case, the projection to $\Sigma$ looks much like it did in Case 2, but now the portions of the $C_i$'s running from $\alpha_0$ to $\alpha_1$ needn't start and end at corresponding points; that is, they twist in the annulus. Inspired by \cite[Lemma 7]{g}, we consider an appropriate power of $\phi$, $\phi^k$. We get $A_\psi$, $\phi(A_\psi)$, \ldots $\phi^{k-1}(A_\psi)$ which patch together to give an annulus $kA_\psi$ between $\alpha_0$ and $\phi^k(\alpha_0)$ of area $k\omega_\phi(A_\psi)$. If we choose $k$ properly, then the $C_i$ running through $kA_\psi$ start and end at corresponding points of $\alpha_0$ and $\phi^k(\alpha_0)$. Thus we're in Case 2 (possibly with an irrelevant Dehn twist in the annulus), and we see that $k\omega_\phi(A_\psi) = 0$, and so $\omega_\phi(A_\psi) = 0$.
\qed

\begin{cor}
\label{iswkmono}
Let $\phi: \Sigma \rightarrow \Sigma$ be in standard form. Then $\phi$ is $\eta$-weakly monotone for all $\eta \in \pi_0(\Gamma(M_\phi))$.\qed
\end{cor}

\begin{cor}
Let $\phi: (\Sigma,\omega)\rightarrow (\Sigma,\omega)$ be a reducible symplectomorphism in standard form in mapping class $h \in \Gamma$. Then $HF_*(\phi) = HF_*(h)$
\end{cor}

\proof
This follows from Corollary \ref{iswkmono} and Corollary \ref{notmonored}.
\qed

\subsection{Floer homology of $\phi_{sm}$ for reducible maps}
\label{flred}

In this section we show how to compute $HF_*(\phi_{sm})$ by showing that $(CF_*(\phi_{sm}), \del)$ splits into a contribution from each component (and annular region between parallel Dehn twists), though what this contribution is may be affected by the behavior of the neighboring regions. The argument follows Gautschi \cite[\S4-5]{g}, with the pseudo-Anosov components and their contribution ``coming along for the ride,'' with the exception of pseudo-Anosov components directly abutting fixed components (cf. Remark \ref{puncspecial} and Section \ref{nielsenred}).

The basic idea is to perturb the fixed components by the Hamiltonian flow of a Morse function (actually a small Morse-Smale function). This Morse function should agree near the boundary with the direction of twisting of a twist component, if any. That is, the Morse function achieves a local maximum at boundaries that twist positively and a local minimum at boundaries that twist negatively. Let $\widehat\phi$ denote the perturbed map. At boundary components of the fixed components where there's a positive twist joining up, a local model for $\widehat\phi$ in a neighborhood of the boundary component is $(q,p) \mapsto (q, p - \delta)$ for $\delta > 0$ small. Similarly, when joining up with a negative twist, the map in a neighborhood of the boundary component is $(q,p) \mapsto (q, p + \delta)$ for $\delta > 0$ small. We let $\Sigma_0$ denote the collection of fixed components, and $\del_{\pm}\Sigma_0$ the collection of components of $\del \Sigma_0$ on which we've joined up with a positive (resp. negative) twist. 

Additionally, the Morse function should agree near the boundary of a pseudo-Anosov component with the Hamiltonian perturbation we chose in Section \ref{unrotpunc}, which we may have to modify slightly by shrinking the size of the perturbation further (i.e. the same process we did already in Section \ref{unrotpunc}, but with a smaller $\epsilon$). We may choose which direction to twist on each side when a pseudo-Anosov component meets a fixed component so long as we do so coherently. The collection of components of the boundary of $\Sigma_0$ which meet a pseudo-Anosov component will be denoted $\del_p \Sigma_0$.

With this setup, we call $\widehat\phi$ a perturbed standard form map. Note that we have left some ambiguity with respect to which direction we're twisting when a fixed component meets a pseudo-Anosov component.

Now Gautschi \cite[Proposition 14]{g} has a neck-stretching argument (actually, he enlarges the symplectic form) to show that different components of $\Sigma \backslash N$ don't interact because flow lines originating at fixed points of Type Ia and Type Ib can't approach the boundary of a component except where a Type IIId boundary abuts a Type Ia component (see the proof of Lemma \ref{cfsplits}, below, for other fixed points). This uses the fact that  $\omega_{\phi}: H_1(\Gamma(M_\phi)) \rightarrow \RR$ is the zero map, i.e. our Lemma \ref{noperiods}. We present a modification of Gautschi's argument in order to show that flow lines between Type Ia and Type IIId fixed points in abutting components don't reach past the Hamiltonian portion of the pseudo-Anosov component corresponding to the Type IIId fixed points. We combine this with an argument of Salamon and Zehnder \cite{sz} to determine these flow lines.

\begin{remark}
Hutchings and Sullivan have an argument involving intersection positivity with the $J$-holomorphic foliation coming from $\RR \times T$ for $T$ a torus in $M_\phi$ which is invariant under the mapping torus flow \cite[Lemma 3.11]{hs}. They consider only such tori coming from fixed circles $C \subset \Sigma$, but the argument extends to collections of circles which are permuted. This can also be used, together with an argument as in \cite[\S3.2-3]{jg} or Case B of the proof of Lemma \ref{noperiods}, to show that flow lines originating at fixed points of Type Ia and Type Ib can't approach the boundary of a component except where a Type IIId boundary abuts a Type Ia component. Their argument seems to be morally equivalent to Gautschi's, with both essentially getting at a local energy contribution. We don't have an argument along the lines of Hutchings and Sullivan's for Lemma \ref{neck}, below, but one may be possible.
\end{remark}

\begin{lemma}
\label{neck}
Consider the subset $\Sigma_0'$ of $\Sigma$ consisting of Type Ia regions which abut Type IIId fixed points together with neighborhoods of the corresponding Type IIId boundaries (including the fixed points) in the appropriate pseudo-Anosov components. Then there exists a modification of the Hamiltonian on the interior of $\Sigma_0'$ such that the differentials between fixed points in $\Sigma_0'$ are given by gradient flow lines of the Hamiltonian.
\end{lemma}

\proof
We perform the same analysis at each Type IIId boundary.

We assume the setup of Step 4 of Section \ref{unrotpunc}. Our plan is to modify what is done in that step. We multiply by $\epsilon_1$ for $r < R$. For $\frac{5R}{4} < r < \frac{3R}{2}$, we multiply by $\epsilon_2$, and we multiply by $1$ for $r > 2R$. We assume that $0 \leq \epsilon_i \leq 1$. We furthermore multiply our Hamiltonian on $\Sigma_0$ by $\epsilon_1$. The plan is to obtain a local energy contribution whenever a differential passes through the region $\frac{5R}{4} < r < \frac{3R}{2}$ and then to use Gautschi's trick of varying the symplectic form to increase the local energy contribution without bound while keeping constant the energy of differentials between fixed points in $\Sigma_0'$.

\emph{Part 1}: We derive bounds on $|\nabla H|$. Fix $R$. We use the standard euclidean metric. We perform a brief calculation to bound $|\nabla H|$ from above and below for $R < r < 2R$. Suppose $H$ is of the form $\kappa \mu r^2 \cos(p\theta)$ (for some $\kappa$) in a subset of this region. Then $$|\nabla H| = \kappa \mu \sqrt{4r^2\cos^2(p\theta) + p^2 r^6 \sin^2(p\theta)}.$$ Using $$\frac{1}{2} \leq \textrm{min}(\cos^2(p\theta),\sin^2(p\theta)) \leq \textrm{max}(\cos^2(p\theta),\sin^2(p\theta)) \leq 1,$$ we have:

\beqa
\kappa\mu\sqrt{\frac{1}{2}\textrm{min}(4r^2,p^2r^6)} &\leq |\nabla H| \leq& \kappa\mu\sqrt{4r^2 + p^2r^6} \\
\Longrightarrow \kappa\mu\frac{1}{\sqrt{2}} \min(2r,p r^3) & \leq |\nabla H| \leq & \kappa\mu (2r + pr^3) \\
\Longrightarrow \kappa\mu\frac{1}{\sqrt{2}}\min(2R,p R^3) & \leq |\nabla H| \leq & \kappa\mu (4R + 8pR^3)
\eeqa

We must also bound the derivative in the region $R < r < \frac{5R}{4}$. When $r = \frac{5R}{4}$, we're supposed to have $H(r) = \epsilon_2\mu r^2\cos(p\theta)$. No matter what $\epsilon_1$ is, we'll have to increase $H$ by a quantity of at most $\epsilon_2\mu \frac{25}{16} R^2$ over an interval of length $\frac{R}{4}$. We can do this keeping $|\nabla H|$ less than twice the ratio, namely $|\nabla H| \leq \epsilon_2\mu\frac{25}{4}R$.

\emph{Part 2}: We seek the local energy contribution. Consider our differential $u(s,t)$ and suppose that $u(s,t)$ achieves the $r$-value $\frac{3R}{2}$. We obtain an energy estimate for such $u$ and then use Gautschi's trick of varying the symplectic form in a neighborhood of the annulus $\frac{5R}{4} < r < \frac{3R}{2}$ to show that a differential between two fixed points in $\Sigma_0'$ cannot reach beyond $\frac{3R}{2}$ at any point.

For $s$ in the interval between $s_1$ and $s_2$, we have that $u(s,0)$ is between $r$-values $\frac{5R}{4} + \frac{R}{16}$ and $\frac{3R}{2} - \frac{R}{16}$. Thus we have bounds of the form $\epsilon_2 K_1 \geq |\nabla H| \geq \epsilon_2 K_2$ (where $K_1$ and $K_2$ are some constants, because $\mu$ and $R$ are fixed).

To find a lower bound on $d(\phi(x),x)$ for $x \in \Sigma$ with an $r$-value of between $\frac{5R}{4} + \frac{R}{16}$ and $\frac{3R}{2} - \frac{R}{16}$, we bound the second derivative of $H$ (that is, the Hessian) from above. It is not hard to see that this is bounded in terms of $\epsilon_2$, $\mu$, and $R$, and the bound in terms of $\epsilon_2$ is linear. Thus we have a bound of the form $|\textrm{Hess}(H)| \leq \epsilon_2 L$. 

The maximum displacement of $\phi$ at time $1$ is bounded by $\epsilon_2 K_1$ (the upper bound on $|\nabla H |$). Thus the maximum change in $\nabla H$ is $(\epsilon_2)^2 K_1 L$. Thus for sufficiently small $\epsilon_2$, this is smaller than $\frac{1}{2}\epsilon_2 K_2$. Thus the displacement of the flow at time $1$ is at least $\frac{1}{2}\epsilon_2 K_2$ (by the above lower bound on $|\nabla H|$).

Thus $$\int_0^1 |\del_t u(s,t)| dt \geq \frac{1}{2}\epsilon_2 K_2$$ for $s$ between $s_1$ and $s_2$. By H\"older's inequality, we have $$\int_0^1  |\del_t u(s,t)|^2 dt \geq \frac{1}{4}(\epsilon_2)^2 (K_2)^2$$ for $s$ between $s_1$ and $s_2$. Thus $$\int_{-\infty}^\infty\int_0^1 |\del_t u(s,t)|^2 dt\ ds \geq \frac{1}{4}(\epsilon_2)^2 (K_2)^2 |s_1 - s_2|.$$

We now use Gautschi's trick. If we multiply the symplectic form and the metric by $M >> 0$ (while keeping the almost complex structure $J$ constant) on the entire pseudo-Anosov component minus a region near the boundary. On the region between the compact pieces of the level set through the $p$ saddle points and the boundary, we multiply the symplectic form and metric by a function $M(s)$ constant on level sets $H = s$ which quickly reaches one (and similarly at any other boundary component). In particular, we've multiplied the symplectic form and metric by $M$ in a neighborhood of the annulus $\frac{5R}{4} < r < \frac{3R}{2}$. The holomorphic curves all remain exactly the same because we haven't changed the complex structure, but we get a better energy estimate in terms of the new metric, with a factor of $M$ showing up:

$$E(u) = \int_{-\infty}^\infty\int_0^1 |\del_t u(s,t)|^2 dt\ ds \geq \frac{1}{4}M(\epsilon_2)^2 (K_2)^2 |s_1 - s_2|.$$

In addition, by Lemma \ref{noperiods}, our action functional is exact and so the energy of a holomorphic curve is given by the difference in the action between its two endpoints. If we temper $M(s)$ from $M$ to $1$ quickly enough, we can make the difference in action between any two fixed points in $\Sigma_0'$ arbitrarily close to the difference in action before we expanded the symplectic form and metric because the region in which we have modified the symplectic form does not separate the fixed points from each other in the sense that they can be connected by a topological cylinder which avoids the region. The action of this topological cylinder is then nearly unchanged, and is equal to the energy of any holomorphic curve $u$ between the two fixed points. Thus as we vary $M$, $E(u)$ must remain (nearly) constant, contradicting the above inequality. Thus differentials between fixed points in $\Sigma_0'$ do not go past $r = \frac{3R}{2}$.

\emph{Part 3}: Now we use the argument of Salamon and Zehnder \cite[Theorem 7.3]{sz} to show that for $\epsilon$ small enough, the only differentials in the region $\Sigma_0$ together with everything within radius $\frac{3R}{2}$ (we call this combined region $\Sigma_0'$) are those corresponding to gradient flow lines of our Hamiltonian (which we make Morse-Smale\footnote{Morse-Smale functions are dense in the space of Morse functions, and so we can do this without disturbing any of our inequalities or critical points qualitatively.}).

The main point is that making $\epsilon_1$ and $\epsilon_2$ small makes $|\nabla H|$ small in the region goverened by $\epsilon_1$ and $\epsilon_2$ --- that is, $\Sigma_0'$ (recall that $|\nabla H| \leq \epsilon_2\mu \frac{25}{4} R$ in the region $R\leq r \leq \frac{5R}{4}$, so $|\nabla H|$  is goverened linearly in this entire region by $\epsilon_1$ and $\epsilon_2$). We note that on $\Sigma_0'$ we may assume that $\nabla H$ is determined up to a scalar multiple (namely, $\epsilon_2$) by keeping the ratio between $\epsilon_1$ and $\epsilon_2$ constant.

We briefly recap Salamon and Zehnder's argument in order to show that it applies in our situation. They show linearization of the $\overline\del$ operator is onto for every differential which corresponds to a gradient flow line \cite[Theorem 7.3 (1)]{sz} for $\nabla H$ sufficiently small. This is a local condition and thus carries over to our setting.

They then argue as follows: if there's a differential for $\tau H$ for $\tau > 0$ arbitrarily close to zero, we can take a sequence of these with $\tau \rightarrow 0$. The energy of such goes like $\tau E$, where $E$ is the maximum energy possibly when using $H$.

By Gromov compactness we get a sequence converging to a broken differential. In our setting, we should add that this broken differential must stay in $\Sigma_0'$ by the reasoning in Part 2 and the fact that there are no fixed points in the same Nielsen class as fixed points in $\Sigma_0'$ which are outside $\Sigma_0'$. The rest of the argument can now be repeated verbatim.
\qed

\begin{lemma}[{cf. \cite[Lemma 16]{g}}]
\label{cfsplits}
For any generic path of complex structures $J_t$ on $\Sigma$, the Floer homology chain complex $(CF_*(\phi), \del_{J_t})$ splits into a sum of chain complexes $(\mathcal{C}_i,\del_i)$ for each component of $\Sigma\backslash N$. These chain complexes are the Floer homology chain complex of the component $S$ with boundary, except that when a pseudo-Anosov component abuts a fixed component, we include the boundary fixed points as part of the fixed component and not as part of the pseudo-Anosov component.

\end{lemma}

\proof
Fixed points in different components of $\Sigma\backslash N$ are in different Nielsen classes by Section \ref{nielsenred} (except where a Type Ia region meets a Type IIId boundary). The moduli spaces of flow lines between different Nielsen classes are empty, as discussed in the proof of Theorem \ref{nodiff}, so indeed the chain complex splits. The neck stretching argument of Gautschi \cite[Proposition 14]{g} mentioned earlier implies that the flow lines don't reach the boundary of the component for components of Type Ia and Type Ib (unless a Type Ia region meets a Type IIId boundary). For all other components, there are no possible differentials because fixed points of Type IIa, Type IIb, Type IIIa, and Type IIIc are all the only fixed point in their Nielsen class by Section \ref{nielsenred}, and fixed points of Type IIIb and fixed points of Type IIId not abutting a component of Type Ia are in a Nielsen class with only points of the same index by Sections \ref{unrot}, \ref{unrotpunc}, and \ref{nielsenred}. For Type Ia components abutting Type IIId boundaries, we use Lemma \ref{neck}.
\qed

This splitting together with the fact that we've identified the summands with the Floer homology of the components $S$ means that we can cite the rest of Gautschi's calculation and add on the contributions from pseudo-Anosov components and from Type Ia components which abut Type IIId fixed points.

First, some notation. Recall the notation of $\Sigma_0$ for the collection of fixed components as well as the three types of boundary $\del_{\pm} \Sigma_0$ and $\del_p \Sigma_0$. Additionally let $\Sigma_1$ be the collection of periodic components and let $\Sigma_2$ the the collection of pseudo-Anosov components with punctures (i.e. before any perturbation in Section \ref{unrotpunc}) instead of boundary components wherever there is a boundary component that meets a fixed component.

We further subdivide $\Sigma_0$. Let $\Sigma_a$ be the collection of fixed components which don't meet any pseudo-Anosov components. Let $\Sigma_{b,p}$ be the collection of fixed components which meet one pseudo-Anosov component at a boundary with $p$ prongs. In this case, we assign the boundary components to $\del_- \Sigma_0$ (this is an arbitrary choice).
Let $\Sigma_{b,p}^\circ$ be the collection of the $\Sigma_{b,p}$ with each component punctured once. Let $\Sigma_{c,q}$ be the collection of fixed components which meet at least two pseudo-Anosov components such that the total number of prongs over all the boundaries is $q$. In this case, we assign at least one boundary component to $\del_+ \Sigma_0$ and at least one to $\del_- \Sigma_0$ (and beyond that, it does not matter).

\begin{thm}
\label{redthm}
Let $\widehat\phi$ be a perturbed standard form map $\Sigma\rightarrow \Sigma$ in a reducible mapping class $h \in \Gamma$ with choices of the signs of components of $\del_p\Sigma_0$ and notation as in the preceding paragraphs. Then
\beqa HF_*(h) = HF_*(\widehat\phi) &\cong& H_{*(\mathrm{mod}\,2)}(\Sigma_a,\del_- \Sigma_a;\ZZ/2) \\
&& \oplus \bigoplus_p \left(H_{*(\mathrm{mod}\,2)}(\Sigma_{b,p}^\circ, \del_- \Sigma_{b,p}; \ZZ/2) \oplus (\ZZ/2)^{(p-1)\left|\pi_0(\Sigma_{b,p})\right|}\right) \\
&& \oplus \bigoplus_q \left(H_{*(\mathrm{mod}\,2)}(\Sigma_{c,q}, \del_- \Sigma_{c,q};\ZZ/2) \oplus (\ZZ/2)^{q\left|\pi_0(\Sigma_{c,q})\right|} \right) \\
&& \oplus\  \ZZ/2^{\Lambda(\widehat\phi|_{\Sigma_1})} \oplus CF_*(\widehat\phi|_{\Sigma_2}),\eeqa
where $\Lambda(\widehat\phi|_{\Sigma_1})$ is the Lefschetz number of $\widehat\phi|_{\Sigma_1}$, the $\ZZ/2^{\Lambda(\widehat\phi|_{\Sigma_1})}$ summand is all in even degree, the other two $\ZZ/2$ summands are all in odd degree, and $CF_*(\widehat\phi|_{\Sigma_2})$ denotes the chain complex for $\widehat\phi$ on the components $\Sigma_2$ (where, as mentioned above, we've left off the Type IIId fixed points which abut a fixed component).

Furthermore, $H_1(\Sigma)$ acts as zero except on the first summand, where it acts by the intersection product. That is, given a component $S$ of $\Sigma_a$ with $\del_+ S = 0$, the action of $H_1(\Sigma)$ is by intersection products for its image under the map to $H_1(\Sigma,\Sigma\setminus S) = H_1(S,\del_- S)$, and given a component $S$ with $\del_- S = 0$, the action of $H_1(\Sigma)$ is by intersection products with the inclusion of $H_1(S)$ to $H_1(\Sigma)$, noting that this action preserves $H_1(S)$.
\end{thm}

\proof
The first summand and the $\ZZ/2^{\Lambda(\widehat\phi|_{\Sigma_1})}$ summand are as in \cite[Theorem 1]{g}. The last summand comes from the fact that there are no differentials in the Floer chain complex on the pseudo-Anosov components as in the proof of Theorem \ref{nodiff}. (We've left off the Type IIId fixed points which meet fixed components and will deal with them elsewhere.)

The sum over $q$ arises in the same manner as the first summand with the extra note that if we choose a Morse function on $\Sigma_c$ which has no maxima or minima, which we can do because every component of $\Sigma_c$ has both a positive and a negative boundary component, then there are no differentials between the fixed points in each of the $\Sigma_c$ components and the various Type IIId fixed points (which correspond to the second term in the sum over $q$).

The sum over $p$ arises in this same manner, except that we choose a Morse function on each component of $\Sigma_b$ with one maximum and no minima (which we can do due to its one positive boundary component). Then one of the Type IIId fixed points cancels with the maximum (they form a cancelling pair), resulting in the homology of the punctured component and in one fewer fixed point in the second half of the sum over $p$.

For the action of $H_1(\Sigma)$, we use the extrinsic or algebraic formulation given in Section \ref{product} for $HF_*(\widehat \phi)$, using Proposition \ref{productinvariant} to see that our calculation holds for $HF_*(h)$. The action preserves Nielsen classes because it is defined in terms of flow-lines so it acts on the summands in the statement of the theorem. It acts as zero on the last two because there are no flow-lines between generators of those summands. Furthermore, it can only act nontrivially on a summand which is nonzero in both degrees mod two, which now leaves only the first. There the correspondence of Morse flow-lines and Floer flow-lines shows that the action is the same as the action of $H_1(\Sigma)$ on Morse homology, which is the same as the intersection product.
\qed

\begin{remark}
In comparing this result with \cite{s1} or \cite{g}, note that we use symplectic Floer homology instead of symplectic Floer cohomology.
\end{remark}

\section{Train tracks and combinatorially computing $HF_*(h)$}
\label{tt}

In order to compute $HF_*(h)$ for a pseudo-Anosov or reducible mapping class from Theorem \ref{nodiff} or Theorem \ref{redthm}, we need to understand the canonical singular representative $\phi_{sing}$ of the pseudo-Anosov mapping class, or of the pseudo-Anosov pieces of the reducible map. In particular, we must understand its fixed points, its singular points, and how each singular fixed point is mapped to itself (i.e. whether it is rotated). All this data is combinatorially computable from the action of a representative of the mapping class on an invariant train track. It is the purpose of this section to describe how this computation works. We note that by work of Bestvina and Handel \cite{bh}, one can algorithmically produce an invariant train track from a representation of $\phi_{sing}$ as a collection of Dehn twists (or, if $\phi_{sing}$ is reducible, produce reducing curves from the Dehn twists, and then repeat with the pseudo-Anosov pieces). Additionally, there is a large collection of examples due to Penner \cite{p} for which producing the invariant train track from the Dehn twists is trivial. We review Penner's construction and compute some examples based on it in the next section. In this section as well, we follow Penner \cite{p}. See also \cite{bh}, \cite[Epilogue]{ph}, \cite[\S9]{m}, \cite{flp}, \cite{t}.

A \emph{train track} $\tau$ on a surface $\Sigma$ is a collection of simple closed curves and graphs in which every vertex has valence at least three embedded in the interior of $\Sigma$. It is required to be a smooth embedding on each closed edge (and curve), both of which we call a \emph{branch}, and at the vertices, called \emph{switches}, the tangent vectors for each edge must be collinear and in addition they must not be all of the same sign (that is, the track must not ``stop''). The train track separates the surface into various faces, which are components of $\Sigma\backslash \tau$. These faces have straight vertices (which we do not count as vertices of the face) and cusp vertices (which we do). We require each face to have negative Euler characteristic, with cusps counting as $-\frac{1}{2}$. Note that the sum of the Euler characteristics of the faces (with this cusp contribution) equals the Euler characteristic of the surface. 
We say that a track $\tau$ \emph{fills} $\Sigma$ if every component of $\Sigma\backslash\tau$ is either a disk or a boundary parallel annulus (with some number of cusps).

For a train track $\tau$, define the \emph{branch space} $V_\tau$ to be the vector space generated by the branches $b_i$ of $\tau$. We say that a track $\tau$ \emph{carries} a track $\tau'$ if we can isotop $\tau'$ such that it lies entirely inside $\tau$. We may assume that switches of $\tau'$ are taken to switches of $\tau$ (but not every switch of $\tau$ will have a switch of $\tau'$ on it, even if $\tau'$ passes through that switch). In this case, we get an \emph{incidence matrix} $M: V_\tau \rightarrow V_{\tau'}$ given by $M b_i = \sum_j M_{ij} c_j$ for branches $b_i$ of $\tau$ and $c_j$ of $\tau'$, where $M_{ij}$ is the count of ``how many times $c_j$ passes over $b_i$.''

A track $\tau$ is \emph{recurrent} if for each branch $b$ of $\tau$, there is a curve $\gamma$ which is carried by $\tau$ and which runs over $b$. A track $\tau$ is \emph{transversely recurrent} if for each branch $b$ of $\tau$, there is a simple closed curve $\gamma$ which intersects $b$ nontrivially and creates no bigons with $\tau$ (a bigon with cusps on any of its edges doesn't count). If both of these are true, we say the track $\tau$ is \emph{birecurrent}.

Following \cite[Epilogue]{ph}, we consider a map $\phi: \Sigma \rightarrow \Sigma$ and a birecurrent train track $\tau$ such that $\phi(\tau)$ is carried by $\tau$ (this is what is meant by an \emph{invariant train track} for $\phi$). In this case we get an incidence matrix $M: V_\tau \rightarrow V_{\phi(\tau)} = V_\tau$ (where the equality is given by identifying $\phi(b)$ with $b$). For $\phi$ to be pseudo-Anosov, it will be necessary that $\tau$ fills $\Sigma$ and also that (possibly after modifying the track by moves given in \cite{bh}, such as collapsing invariant trees) $M$ is \emph{Perron-Frobenius}, which means that some (positive) power of $M$ has all positive entries. 
It then follows that $M$ has a unique (up to scaling) eigenvector $v = (v_i)$ with all positive entries. Furthermore, its (positive real) eigenvalue $\lambda > 1$ is the eigenvalue of $M$ which is maximum in absolute value. $M^T$ is also Perron-Frobenius, and we let $w = (w_i)$ denote its unique (up to scaling) eigenvector with all entries positive. Its eigenvalue must also be $\lambda$.

From this data we can construct the measured foliations for $\phi_{sing}$. We fatten each branch and let its width be $v_i$ and its length be $w_i$. These automatically satisfy the \emph{switch condition} that the sum of the widths of branches abutting a switch from one direction equals the sum of the widths of the branches abutting the switch from the opposite direction. Additionally, they satisfy a dual switch condition on the $w$'s (see Penner \cite{p}). It then turns out that we can ``zip up'' the fattened track, collapsing the faces and getting a singular point for each face \cite{p}. 

\subsection*{Computing the number of fixed points}

Now we come to the crux of the matter, which is how to compute the number of fixed points based on this description. To warm up we consider a toy model: suppose $R$ is a rectangle in the plane (parallel to the $x$- and $y$-axes) and we apply an affine transformation $\Phi$ which stretches by a factor of $\lambda > 1$ in the $x$-direction and shrinks by a factor of $\lambda$ in the $y$-direction. Suppose furthermore that the vertical extent of $\Phi(R)$ is contained within the vertical extent of $\Phi$ and the horizontal extent of $R$ is contained within the horizontal extent of $\Phi(R)$, where in both cases we allow their boundaries to coincide. 
Then $\Phi$ has precisely one fixed point in $R$, as can be readily deduced from the contraction mapping principle in the $x$ and $y$ coordinates separately. Furthermore, if their boundaries coincide either horizontally or vertically, then the fixed point will be on this boundary (and if two boundaries coincide, the fixed point will be at the common vertex). Note also that we again get a unique fixed point if the affine transformation stretches in the $x$-direction, contracts in the $y$-direction, and is then rotated by $180^\circ$ (i.e. is negative hyperbolic instead of positive hyperbolic), subject to the the same conditions on the vertical/horizontal extents.

What this means is that there is one fixed point of $\phi$ every time the branch $\phi(b_i)$ runs over the branch $b_i$. This count is the trace of $M$. Unfortunately, we have overcounted, counted some singular fixed points together with the nonsingular fixed points, and additionally undercounted a specific type of nonsingular fixed point, so we have to make some corrections.

In order to deal with the overcounting and undercounting, we formalize two further pieces of data. These are the tangential and transverse lists. Arbitrarily orient each branch $b_i$. We let $\sigma_i$ denote the (tangential) ordered list given by the $b_j$'s which $\phi(b_i)$ runs over (we also record whether the orientations agree when $\phi(b_i)$ runs over $b_j$). Note that the number of times $b_j$ appears in $\sigma_i$ equals $M_{ij}$. Also, let $\nu_i$ denote the (transverse) ordered list given by the $b_j$'s for which $\phi(b_j)$ runs over $b_i$, in order from left to right (with $b_i$ pointing up; again we record whether their orientations agree). Note that the number of times $b_j$ appears in $\nu_i$ is given by $M_{ji}$.

Let $u_p$ denote the number of unrotated fixed singular points with $p$ prongs and let $r_p$ denote the number of rotated fixed singular points with $p$ prongs. Both of these can be read off from the track $\tau$ together with the lists $\sigma_i$ and $\nu_i$. Rotated fixed singular points do not contribute any to $\tr(M)$, but unrotated fixed singular points contribute $p$ to $\tr(M)$, corresponding to $p$ points which are identified when the fattened track is zipped up (ignoring under- and over-counting that will be dealt with below). 
Note that there can be no other fixed points that are on the boundary of the fattened track (and which would then be identified when this boundary is zipped up) because the track is stretched by $\lambda$.

Next we have overcounting of nonsingular fixed points when there is a fixed point at the end of a (fattened) branch, i.e. living over a switch, for then it is counted twice, once in each of the two branches it lies in. 
Let $s_i$ equal one-half the number of times that $b_i$ is the first entry in $\sigma_i$ plus one-half the number of times it's the last (in both cases we require that the orientation of $\phi(b_i)$ agree with the orientation of $b_i$). The number $\sum_i s_i$ represents how many of these sorts of overcounting errors there are.

Finally we have undercounting of nonsingular fixed points when there is a ``flip.'' That is when $b_i$ and $b_j$ share a common switch and, supposing without loss of generality that this is the initial switch for each of $b_i$ and $b_j$, the list $\sigma_i$ begins with $b_j$ and $\sigma_j$ begins with $b_i$ (both with compatible orientations). Let $f$ denote the number of flips. 

\begin{thm}
\label{pacomb}
For $\phi$ pseudo-Anosov on a closed surface $\Sigma$, the rank of $HF_*(\phi)$ in the above setup is $tr(M) + f - \sum_i s_i + \sum_p (r_p - u_p)$.
\end{thm}

\proof
By Theorem \ref{nodiff}, we get a contribution of $\sum_p (p-1) u_p$ from the unrotated singular points, a contribution of $\sum_p r_p$ from the rotated singular points, and the above discussion shows that the number of normal fixed points is $\tr(M) + f - \sum_i s_i - \sum_p p\ u_p$. Adding together these three quantities, we get the desired result.
\qed

We could get the number of even generators and the number of odd generators from the lists $\sigma_i$ and $\nu_i$ with orientations attached, if desired. Alternatively, we can simply use the Lefschetz fixed point theorem to get the difference between the number of odd and the number of even generators, and combine that with the above theorem. 

\section{Penner's construction and computation of examples}
\label{exs}

For a simple closed curve $C$ on $\Sigma$, let $\tau_{C}$ be the positive Dehn twist about $C$. Consider a collection $\{C_i\}$ of disjoint embedded curves on $\Sigma$, and $\{D_j\}$ another such collection. Let $\phi$ be any composition of positive Dehn twists along the $C_i$'s and negative Dehn twists along the $D_j$'s such that every curve listed shows up at least once. Penner \cite{p} shows that that if every region of the complement $\Sigma -\coprod_i C_i - \coprod_j D_j$ is a disk with more than two edges, then the mapping class $[\phi]$ is pseudo-Anosov. Furthermore, if every curve listed shows up at least twice, then we have an explicit invariant train track for $\phi$ whose incidence matrix $M$ is Perron-Frobenius.\footnote{This train track will still be invariant if some curves show up only once, but the incidence matrix might not be Perron-Frobenius and we might need to perform some moves as in \cite{bh}.}

This invariant train track $\tau$ is simply take the union of all the curves $C_i$ and $D_j$ with their intersection points modified to be switches: at an intersection point of $C_i$ with $D_j$, have the switch curve to the right as viewed along $C_i$ (and thus to the left as viewed along $D_j$). This convention guarantees that the train track is invariant under positive Dehn twists along the $C_i$'s and under negative Dehn twists along the $D_j$'s.

At an intersection point of $C_i$ with $D_j$, there are two ends of $D_j$ meeting $C_i$, and these (possibly) correspond to separate branches of the train track $\tau$. If we perform a Dehn twist along a curve just to one side of $C_i$, then one of these branches is modified to one which goes all the way around $C_i$, and if we Dehn twist on the other side of $C_i$, the other is. We will be twisting along $C_i$ at least twice, so we may ensure that we twist on each side of $C_i$ at least once. This will ensure that the incidence matrix $M$ is Perron-Frobenius.

The incidence matrix $M$ is simply the product of the incidence matrices for each curve we twist along (the matrix depends on the side of the curve we twist on). Thus $\tr(M)$ can be computed as the trace of this product of matrices. The other terms in Theorem \ref{pacomb} are also combinatorial. We have:

\begin{thm}
\label{pennerex}
Let $\{C_i\}$ and $\{D_j\}$ be collections of noncontractible simple closed curves on $\Sigma$ with the $C_i$ mutually disjoint and the $D_j$ mutually disjoint such that every region of the complement $\Sigma-\coprod_i C_i - \coprod_j D_j$ is a disk with more than two edges.

Let $\phi = \prod_j \tau_{D_j}^{-2} \prod_i \tau_{C_i}^2$ and let $h$ be the mapping class of $\phi$. Then $HF_*(\phi) = HF_*(h)$ has rank $$4\sum_{i,j} \left|C_i \cap D_j\right|^2 + 2g - 2.$$
\end{thm}

\proof
We calculate each of the terms in the formula given in Theorem \ref{pacomb}. In order to find $\tr(M)$, we express $M$ as a product of the incidence matrix for $\prod_j \tau_{D_j}^{-2}$ with the incidence matrix for $\prod_i \tau_{C_i}^2$.

First we fix some notation. We have curves $C_i$ and $D_j$ and branches $c_\alpha$ and $d_\beta$ which lie on some $C_i$ and some $D_j$, respectively. We let $C_{[\alpha]}$ denote the curve $C_i$ such that $c_\alpha$ lies on $C_i$, and likewise for $D_{[\beta]}$. We will write our incidence matrices in block form with the entries associated to the $c_\alpha$ and the $d_\beta$ separated.

The incidence matrix for $\prod_i \tau_{C_i}^2$ is of the form
$$M_C = \left[\begin{array}{cc}
I & 0 \\
A &  I
\end{array}\right]$$
where $A$ is a matrix with entries $A_{\beta\alpha}$ equal to the number of times the branch $d_\beta$ passes over $c_\alpha$. This equals the number of times $d_\beta$ passes over the curve $C_{[\alpha]}$, which equals the number of times the branch $d_\beta$ intersects the curve $C_{[\alpha]}$, that is, how many of its ends lie on $C_{[\alpha]}$, with each end counting separately. We denote this as $\#(d_\beta \cap C_{[\alpha]})$, though this may not equal their set-theoretic intersection number when the two ends of $d_\beta$ coincide. We get a contribution from each end of $d_\beta$ because we've twisted on either side of each $C_i$. Likewise, the incidence matrix for $\prod_j \tau_{D_j}^{-2}$ is of the form
$$M_D = \left[\begin{array}{cc}
I & B \\
0 &  I
\end{array}\right]$$
where $B$ is a matrix with entries $B_{\alpha\beta}$ equal to the number of ends of the branch $c_\alpha$ which lie on $D_{[\beta]}$.

We have
$$M = M_D M_C = 
\left[\begin{array}{cc}
I & B \\
0 &  I
\end{array}\right]
\left[\begin{array}{cc}
I & 0 \\
A &  I
\end{array}\right]
= \left[\begin{array}{cc}
I + BA & B \\
A &  I
\end{array}\right].$$ Thus \beqan\label{eq_tr} \tr(M) = b + \sum_{\alpha,\beta} B_{\alpha\beta}A_{\beta\alpha},\eeqan where $b$ is the total number of branches.

We claim that \beqan\label{eq_sum} \sum_{\alpha\beta} B_{\alpha\beta}A_{\beta\alpha} = 4\sum_{i,j} \left|C_i \cap D_j\right|^2.\eeqan
To see this, note that 
\beqa
\sum_{\beta\in[\beta]} B_{\alpha\beta} A_{\beta\alpha} &=& \sum_{\beta\in[\beta]} \#\!\left(c_\alpha \cap D_{[\beta]}\right)\#\!\left(d_\beta \cap C_{[\alpha]}\right) \\
&=& \#\!\left(c_\alpha \cap D_{[\beta]}\right)2\#\!\left(D_{[\beta]} \cap C_{[\alpha]}\right),\eeqa
where the factor of two comes from the fact that each intersection point of $D_{[\beta]}$ with $C_{[\alpha]}$ is associated to two ends of branches in $D_{[\beta]}$. Now summing over $\alpha\in[\alpha]$ gives equation \ref{eq_sum}. We furthermore note that the number of branches $b$ is twice the number of intersections, i.e.~$b = 2\sum_{i,j}\left|C_i \cap D_j\right|$, because each intersection involves four branch ends (of which there are two per branch).

It remains to compute the remaining terms in Theorem \ref{pacomb}. First we consider the start and end of each branch (from its ordered list $\sigma$). Consider the branch $c_\alpha$. It is unchanged by twisting along the $C_i$'s. When we twist along the $D_j$'s, $c_\alpha$ covers the $D$-curve at its initial end, then itself, then the $D$-curve at its final end. In particular, $s_\alpha$ is zero and there are no contributions to flips from the $c_\alpha$'s alone. Now consider the branch $d_\beta$. When we twist along the $C_i$'s, it is modified and covers the $C$-curve at its initial end, then itself, then the $C$-curve at its final end, but also it starts and ends on those two $C$-curves at its initial and final switches, and in the same orientation (at the switch) as itself at the switches. Call these initial and final branches $c_{init}$ and $c_{fin}$. Then when we twist along the $D_j$'s, by the same reasoning, the initial switch that $c_{init}$ runs over is $d_\beta$, and the final switch that $c_{fin}$ runs over is $d_\beta$, each with the appropriate orientation. Thus the ordered list $\sigma_\beta$ begins and ends with $d_\beta$ (with the correct orientations). Thus $s_\beta = 2/2 = 1$ and also there are no flips.

Thus we get a (negative) contribution of $\sum_i s_i = \sum_\beta 1$. This equals the number of branches in the $D$-curves. Because each switch has two $d$-branches and two $c$-branches incident to it, and each branch has two ends, we get that this equals the number of switches, which is $\sum_{i,j}\left|C_i \cap D_j\right|$.

It remains to examine the singularities. These correspond to regions in the complement of the $C$'s and $D$'s. We can determine the number of these using Euler characteristic: the number of switches minus the number of branches plus the number of these regions equals the Euler characteristic of the surface, which is $2-2g$. We know the number of switches is $\sum_{i,j}\left|C_i \cap D_j\right|$ and the number of branches is $2\sum_{i,j}\left|C_i \cap D_j\right|$. Thus there are $2-2g + \sum_{i,j}\left|C_i \cap D_j\right|$ singularities.

We claim that each of these singularities is fixed and unrotated. To determine this requires understanding the transverse list $\nu_\beta$ for each branch $d_\beta$: the initial and final branches in this list (with orientation) determine which regions are mapped to the regions on either side of $d_\beta$, and how rotated they are.\footnote{Note that because each region borders some $d_\beta$'s, we won't need to consider the $c_\alpha$'s separately.} Thus we must show that the initial and final branches in this list are both $d_\beta$ with the correct orientation.

When we perform a Dehn twist along $D_{[\beta]}$, the branches that cover a given branch $d_\beta$ are those $c_\alpha$'s which meet $D_{[\beta]}$. These cover $d_\beta$ on the side corresponding to the side of $D_{[\beta]}$ we've twisted on. The outermost branches covering $d_\beta$ (i.e.~the initial and final branches in $\nu_\beta$) are those which meet $d_\beta$ at a switch which are on the opposite side of the switch from $d_\beta$. Thus the outermost branches covering $d_\beta$ are the outermost branches covering these $c_\alpha$'s after twisting along $C_{[\alpha]}$ (for then these will be carried along with $c_\alpha$ and cover $d_\beta$ when we twist along $D_{[\beta]}$). These are likewise the $d$-branches meeting the $c_\alpha$ at a switch on the opposite side of the switch. One of these is $d_\beta$, and we can see that in fact this branch is the outermost when $c_\alpha$ covers $d_\beta$, as desired.

Summing up, we've shown that 
\beqa
HF_*(\phi) &=& \tr(M) + f - \sum_i s_i + \sum_p(r_p - u_p) \\
&=& \left(2\sum_{i,j}\left|C_i \cap D_j\right|+ 4\sum_{i,j} \left|C_i \cap D_j\right|^2\right) + 0 \\
&& - \sum_{i,j}\left|C_i \cap D_j\right| - \left(2-2g + \sum_{i,j}\left|C_i \cap D_j\right|\right) \\
&=& 4\sum_{i,j} \left|C_i \cap D_j\right|^2 + 2g - 2
\eeqa

\qed

\begin{remark}
This easily extends to $\prod_j \tau_{D_j}^{-2\ell} \prod_i \tau_{C_i}^{2k}$. Only $\tr(M)$ is changed (we replace $A$ by $kA$ and $B$ by $\ell B$). We get $HF_*(\phi) = 4k\ell\sum_{i,j} \left|C_i \cap D_j\right|^2 + 2g - 2$. In fact, any composition of positive powers of $\prod_j \tau_{D_j}^{-2}$ and positive powers of $\prod_i \tau_{C_i}^{2}$ can be computed in this manner. The correction terms again remain the same, so we only need to know how to compute $tr(M)$. 
\end{remark}

\end{document}